\documentclass[11pt]{amsart}

\usepackage{amsxtra}
\usepackage{amssymb}

\newcommand\la{\langle}
\newcommand\ra{\rangle}

\newcommand\dd{{\mathfrak d}}
\renewcommand\aa{{\mathfrak a}}
\newcommand\bb{{\mathfrak b}}

\newcommand\uu{{\mathfrak u}}
\newcommand\hh{{\mathfrak h}}
\newcommand\rr{{\mathfrak r}}
\newcommand\nn{{\mathfrak n}}
\newcommand\ggo{{\mathfrak g}}
\newcommand\aff{\mathfrak {aff}}

\newcommand\vv{{\mathfrak v}}

\newcommand\CC{\mathbb C}

\newcommand\RR{\mathbb R}

\newcommand\Der{\operatorname{Der}}
\newcommand\ad{\operatorname{ad}}

\newcommand\tr{\operatorname{tr}}

\theoremstyle{plain}

\newtheorem{thm}{Theorem}[section]
\newtheorem{lem}[thm]{Lemma}
\newtheorem{prop}[thm]{Proposition}
\newtheorem{cor}[thm]{Corollary}

\theoremstyle{definition}

\theoremstyle{remark}
\newtheorem*{remark}{Remark}
\newtheorem*{remarks}{Remarks}

\title{Product structures on four dimensional solvable Lie algebras}

\author{A.~Andrada}

\author{M.~L.~Barberis}

\author{I.~G.~Dotti}

\author{G.~P.~Ovando}

\address{CIEM, FaMAF, Universidad Nacional de C\'ordoba, Ciudad Universitaria,
(5000) C\'ordoba, Argentina}
\thanks{The authors were partially supported
by CONICET and SECYT-UNC (Argentina).}


\begin{document}

\begin{abstract} It is the aim of this work to study product structures on four
dimensional solvable Lie algebras. We determine all possible paracomplex structures and consider the case when one of the subalgebras is an ideal. These results are applied to
the case of Manin triples and complex product structures. We also analyze the three dimensional subalgebras.
\end{abstract}

\maketitle

\section*{Introduction}

A product structure on a smooth manifold $M$ is an endomorphism $E$ of its tangent bundle  satisfying $E^2=\,$Id together with \begin{equation}
E[X,Y]=[EX,Y]+[X,EY]-E[EX,EY] \quad \text{for all vector fields }
X,Y \text{ on } M.  \end{equation}
A product structure on $M$ gives rise to a splitting of the tangent bundle $TM$
into the Whitney sum of two subbundles $T^{\pm}M$ corresponding to the
$\pm 1$ eigenspaces of $E$. The distributions on $M$ defined by $T^+M$ and $T^-M$
 are completely integrable. When  $T^+M$ and $T^-M$ have the same rank the product
structure is called a paracomplex structure.

Product structures on manifolds were considered by many authors from different points of view. Examples of Riemannian almost product structures were given in \cite{Miq} and a survey on paracomplex geometry can be found in \cite{CFG}. The classification of Riemannian almost product manifolds according to a certain decomposition of the space of tensors was done in
\cite{N}. In \cite{LM} the authors give a new look at singular and non holonomic
Lagrangian systems in the framework of almost product structures.
Complex product structures on Lie groups were considered in \cite{AS}
and \cite{BV}.

In this paper we consider product structures on four dimensional solvable Lie groups.
Such groups provide an important source
of  applications in geometry.  Invariant structures on the group, for instance, special metrics \cite{Al}, \cite{B2}, \cite{D-S}, \cite{F1}, \cite{F2}, \cite{J}, complex and K\"ahler structures \cite{ACFM}, \cite{AFGM}, \cite{O1}, \cite{SJ}, \cite{FG}, hypercomplex and hypersymplectic structures \cite{An}, \cite{B1}, can be read off in $\RR^4$, the universal covering group, giving often explicit descriptions of the corresponding
structure.

A left-invariant product structure on a Lie group is determined by its restriction to the corresponding Lie algebra, considered as the tangent space at the identity. A {\em  product structure} on a Lie algebra $\ggo$ is a linear endomorphism $E:\ggo \longrightarrow \ggo$ satisfying $E^2=\,$Id (and not equal to $\pm $Id) and
\begin{equation}
E[x,y]=[Ex,y]+[x,Ey]-E[Ex,Ey] \quad \text{for all } x,y\in \ggo.
\end{equation}
A product structure on $\ggo$ gives rise to a decomposition of $\ggo$ into
\begin{equation}
\ggo= \ggo_{+} \oplus \ggo_{-},\;\;\; E|\ggo_{+}=\text{Id},\; E|\ggo_{-}=-\text{Id},
\end{equation}
where both, $\ggo_{+}$ and $\ggo_{-}$, are Lie subalgebras of $\ggo$. This will be denoted
$\ggo = \ggo_+ \bowtie \ggo_-$, since the structure of $\ggo$ is that of a double Lie algebra (\cite{LW}). In case both $\ggo_+$ and $ \ggo_-$ have the same dimension we say that $\ggo$ carries a paracomplex structure.

The outline of this paper is as follows. In Section 1 we describe all non-isomorphic four dimensional solvable Lie algebras over $\RR$. This was studied by Mubarakzyanov \cite{Mu} and Dozias \cite{D}. We found citations of the theorems obtained by Mubarakzyanov in \cite{Z}, pp.\ 988 and Dozias in \cite{Ve}, pp.\ 180. We include a proof of the classification theorem since it will be frequently used to obtain the results throughout the article. Appendix II contains comparisons with the tables given by the various authors \cite{Mu}, \cite{D}, \cite{SJ}, \cite{O1}, \cite{Z}.

In Section 2 we consider product structures on four dimensional Lie algebras. We determine all four dimensional solvable Lie algebras admitting a paracomplex structure (see Table~\ref{pc}).
Among these, we study the case when one of the subalgebras is an ideal of $\ggo$. We also exhibit decompositions where one of the subalgebras is three dimensional (see Table~\ref{1,3}).

An important subclass of paracomplex structures is given by Manin triples and complex product structures (see Section 3). A paracomplex structure $\ggo =\ggo_+ \bowtie \ggo_-$ is a Manin triple if there exists a non degenerate invariant symmetric bilinear form on $\ggo$ such that $\ggo_{\pm}$ are isotropic subalgebras. It is shown that there is only one non abelian four dimensional solvable Lie algebra giving rise to a Manin triple. On the other hand, given a product structure $E$ and a complex structure $J$ on $\ggo$ such that $JE=-EJ$, $\{J,E\}$ is called a complex product structure on $\ggo$. We determine all four dimensional solvable Lie algebras admitting complex product structures (see Table~\ref{cps}), giving an alterantive proof of a result by Blazi\'c and Vukmirovi\'c (\cite{BV}).

\

\section{Classification of four dimensional solvable Lie algebras}
In this section we obtain the classification of four dimensional
solvable Lie algebras. The proof follows the lines of
\cite{Mi} for the classification of three dimensional solvable Lie
algebras, that is, we obtain the four dimensional solvable Lie
algebras as extensions of the three dimensional unimodular Lie
algebras $\RR^3$, the Heisenberg algebra $ \hh_3$, the Poincar\'e
algebra $\mathfrak e(1,1)$ or the Euclidean algebra $ \mathfrak
e(2)$. Both, \cite{O1} and \cite{SJ}, obtain the four dimensional
solvable Lie algebras as extensions of nilpotent Lie algebras of
dimension at most three. In Appendix I we exhibit matrix
realizations and Appendix II contains comparisons with the tables
given by the various authors \cite{Mu}, \cite{D}, \cite{SJ},
\cite{O1}, \cite{Z}.

\subsection{Algebraic preliminaries}

A Lie algebra $\ggo$ which satisfies the condition $\tr(\ad(x))=
0$ for all $ x \in \ggo$ will be called a {\it unimodular} Lie
algebra. If $\ggo$ is a Lie algebra, then using the Jacobi identity we see that
$\tr(\ad[x,y])=0$ for all $x, y \in \ggo$. Hence, the map $\chi :
\ggo \to \RR$ defined by
\begin{equation} \label{chi}
\chi(x)=  \text{tr}(\ad(x)), \qquad\qquad x\in \ggo,
\end{equation}
is a Lie algebra homomorphism. In particular, its kernel $  \uu =
\ker(\chi)$ is an ideal containing the commutator ideal
$[\ggo,\ggo]$. The ideal $\uu$ will be called the {\it unimodular kernel} of $\ggo$. It is easy to check that $\uu$ itself is
unimodular.

We now introduce some notation that will be used throughout the paper (compare with \cite{GOV}).
\begin{itemize}
\item [$\aff(\RR)$: ]  $[e_1,e_2]=e_2$, the two dimensional
non-abelian Lie algebra of the group of affine motions of the real
line; \item [$\hh_3$: ] $[e_1,e_2] = e_3$, the three-dimensional
Heisenberg algebra; \item[$\mathfrak r_3$: ] $[e_1,e_2]=e_2,  \;
[e_1,e_3]=e_2+e_3$; \item [$\mathfrak r_{3,\lambda}$: ]
$[e_1,e_2]=e_2,  \; [e_1,e_3]= \lambda e_3 $; \item [$\mathfrak
r'_{3,\lambda}$: ] $[e_1,e_2]=\lambda e_2 - e_3,  \; [e_1,e_3]=
e_2 +\lambda e_3$;
\end{itemize}

\begin{remark}
Observe that $\mathfrak r_{3,-1}$ is the Lie algebra $\mathfrak
e(1,1)$ of the group of rigid motions of Minkowski $2$-space,
$\mathfrak r_{3,0}= \RR \times \aff(\RR)$ and $\mathfrak r_{3,1}$
is the Lie algebra of the solvable group which acts simply and
transitively on the real hyperbolic space $\RR H^3$. Also
$\mathfrak r_{3,0}'$ is the Lie algebra $\mathfrak e (2)$ of the
group of rigid motions of Euclidean $2$-space. Other authors
denote $\aff(\RR)$ by $\mathfrak{sol}_2$ and $\mathfrak e(1,1)$ by
$\mathfrak{sol}_3$.
\end{remark}

We recall the classification of solvable Lie algebras of dimension
$\leq 3$. A proof can be found, for example, in \cite{Mi} or \cite{GOV}.

\begin{thm}
Let $\ggo$ be a real solvable Lie algebra, $\dim\ggo\leq 3$. Then
$\ggo$ is isomorphic to one and only one of the following Lie
algebras: $\; \RR$, $\; \RR^2$, $\;\aff(\RR)$, $\;\RR^3$,
$\;\hh_3$, $\;\mathfrak r_3$, $\; \mathfrak r_{3,\lambda}, \;
|\lambda|\leq 1 \;$ and  $\;\mathfrak r'_{3,\lambda}, \; \lambda
\geq 0$. Among these, the unimodular ones are $\; \RR$, $\;
\RR^2$,  $\;\RR^3$, $\;\hh_3$,  $\; \mathfrak r_{3,-1}$, and $\;
\mathfrak r'_{3,0}$.
\end{thm}

The proof of Theorem \ref{classes} in next section is based on the knowledge of the algebra of derivations of solvable unimodular three dimensional Lie algebras. This is the content of the next lemma, whose proof is straightforward.

\begin{lem} \label{der} The algebra of derivations of $\;\mathfrak e(2),\;\mathfrak e(1,1)$ and $\hh_3$ are
\begin{equation}\label{der1}
 \Der\mathfrak e(2) = \left\{
\begin{pmatrix}   0&0&0 \\  c &a & -b   \\
 d & b & a
\end{pmatrix} : \; a, b, c, d\in \RR \right\}\cong \aff (\CC), \end{equation}
with respect to the basis $e_i, \; i=1,2,3,\;$ such that $[e_1,e_2]=e_3, \;\; [e_1,e_3]=-e_2$;
\begin{equation}\label{der2} \Der\mathfrak  e(1,1)=\left\{
\begin{pmatrix} 0&0&0\\ c & a & 0   \\
d & 0 & b  \end{pmatrix} : \; a, b, c, d\in \RR \right\}\cong
 \aff(\RR)\times\aff(\RR),\end{equation}
with respect to the basis $e_i, \; i= 1,2,3\;$ such that
$[e_1,e_2]=e_2, \;\; [e_1,e_3]=-e_3$;
\begin{equation}\label{der3} \Der\mathfrak h_3  =\left\{
\begin{pmatrix} A & \begin{matrix} 0 \\  0 \end{matrix}  \\
\, b \;\; \; c & \tr A
\end{pmatrix} : \; A \in \ggo \mathfrak l(2,\RR), \; b, c\in \RR \right\},\end{equation}
with respect to the basis $e_i, \;\; i=1,2,3\;$ such that $[e_1,e_2]=e_3$.
\end{lem}

\subsection{Classification theorem}
In this section we obtain all four dimensional solvable Lie algebras as semidirect extensions of three dimensional unimodular Lie algebras. The classification theorem is then reduced to the study of the derivations of these three dimensional algebras. The proof will follow the lines of \cite{Mi} for the three dimensional case, but instead of the rational form, we make use of
the Jordan normal form over ${\RR}$.

Given a Lie algebra $\ggo$ and an ideal $\mathfrak v$ of
codimension one in $\ggo$, let $e_0 \in \ggo \setminus
\mathfrak v$. Then we denote
\begin{equation} \label{eq-sd}
\ggo=\RR e_0 \ltimes_{\varphi}\mathfrak v,
\end{equation}
where $\varphi:\RR e_0\to\Der\mathfrak v$ is a linear map such
that $\varphi(e_0)=\ad(e_0)$. Observe that the splitting of the
short exact sequence
\[ 0 \to \vv \to \ggo \to \RR \to 0, \]
is an immediate consequence of the fact that $\RR$ is one dimensional.

The following result proves the desired decomposition, that is, any four dimensional solvable real Lie algebra is a semidirect product of $\RR$ and a three-dimensional unimodular ideal. Thus this proposition is a first step in the classification  (compare with Proposition~2.1 in
\cite{D-S}):

\begin{prop}\label{ideal}
Let $\ggo$ be a four-dimensional solvable real Lie algebra. Then
there is a short exact sequence
\[ 0 \to \vv \to \ggo \to \RR \to 0, \]
where $\vv$ is an ideal of $\ggo$ isomorphic to either
$\RR^3,\;\hh_3,\;\mathfrak e(1,1)$ or $\mathfrak e(2)$, that is,
$\ggo \cong \RR e_0 \ltimes_{\varphi} \mathfrak v$.
\end{prop}

\begin{proof} Consider the  Lie algebra homomorphism $\chi:\ggo\to\RR$ defined in \eqref{chi}. If $\ggo$ is not
unimodular then its unimodular kernel $\uu$ has dimension three,
therefore it is isomorphic to $\RR^3,\,\hh_3,\,\mathfrak e(1,1)$
or $\mathfrak e(2)$ and the proposition follows with $\vv=\uu$.

We assume now that $\ggo$ is unimodular. The commutator ideal
$\ggo'$ is nilpotent and $\dim \ggo' \leq 3$, hence it follows
that $\ggo'$ is isomorphic to $\{ 0 \},\; \RR, \; \RR^2,\;\RR^3$
or $\hh_3$. In the last two cases the proposition follows by
taking $\vv=\ggo'$. If $\ggo'= \{ 0 \}$ then $\ggo$ is abelian so
that $\vv=\RR^3$ is an ideal of $\ggo$.

If $\ggo'$ is isomorphic to $\RR$, $\ggo' =\RR e_3$, then there
exist elements $e_1, e_2$ in $\ggo$ such that $[e_1,e_2] = e_3$.
The set $e_1,e_2,e_3$ is linearly independent since $\ggo$ is
unimodular. Therefore, the Lie subalgebra generated by $e_1, e_2,
e_3$ is an ideal isomorphic to $\hh_3$.

If $\ggo'$ is isomorphic to $\RR^2$ then  either i) there exists
$x$ not in $\ggo'$ such that $\ad(x)_{|{\ggo'}}$ is non singular,
or ii) for all $x\in \ggo$ the transformation $\ad(x)$ is
singular. Making use of the Jordan form of the corresponding
complex transformation we get in both cases i) and ii), that
$\chi(x) = \lambda _1 + \lambda_2=0$, for $\lambda _i \in \CC$,
$i=1$ or 2. Thus in  case i)  there is a basis of $\ggo'$ such
that the action of $x$ is given as  follows (up to a nonzero
multiple):
\[ \text{a)}\,
\ad(x){|_{\ggo'}}=\left(\begin{matrix} 1 & 0 \\ 0 & -1
\end{matrix}\right)  \qquad  \text{or} \qquad
\text{b)}\,\ad(x){|_{\ggo'}}=\left(\begin{matrix} 0 & 1 \\-1 & 0
\end{matrix} \right)\]
where  case b) corresponds to the eigenvalues $i,-i$. Thus $\RR
x \oplus \ggo'$ is an ideal of $\ggo$ isomorphic to $\mathfrak e(1,1)$ or
to $\mathfrak e(2)$, respectively.

In  case ii), since $\lambda_1$ or $\lambda_2$ is zero, then the
unimodular condition imposes that both eigenvalues vanish and so,
for a fixed $x$ not in $\ggo'$,  there is a basis of $\ggo'$ such
that the action of $\ad(x){|_{\ggo'}}$ takes one of the following forms: \[ \text{a)}\,
\ad(x){|_{\ggo'}}=0  \qquad \text{or} \qquad
\text{b)}\,\ad(x){|_{\ggo'}}=\left(\begin{matrix} 0 & 1\\0 & 0
\end{matrix} \right)\]
Therefore, $\RR x \oplus \ggo'$ is an ideal of $\ggo$
isomorphic to $\RR^3$ in  case a) or  $\hh_3$ in  case b).
This completes the proof.
\end{proof}

\smallskip

The following lemma will be used in the proof of the classification theorem.

\begin{lem} \label{R3}
Let $\ggo_1=\RR e_0\ltimes_{\varphi_1}\RR^3$ and $\ggo_2=\RR e_0
\ltimes_{\varphi_2}\RR^3$ such that $[\ggo_i,\ggo_i]=\RR^3,\;
i=1,2$. Then $\ggo_1 \cong \ggo _2$ if and only if there exists
$\gamma \neq 0$ such that $\varphi_1 (e_0)$ and $\gamma
\varphi_2(e_0)$ are conjugate in GL$(3,\RR)$.
\end{lem}

\begin{proof} Assume first that there exists a Lie algebra isomorphism
$\psi:\ggo_1\to\ggo_2$; then $\psi:\RR^3 \to \RR^3$ and
$\psi(e_0)=\gamma e_0+w$, where $\gamma\in\RR\smallsetminus \{0\}$
and $w\in \RR^3$. If $v\in\RR^3$, we calculate
\begin{align*}
[\psi(e_0),\psi (v)] & = \gamma \varphi _2(e _0) \psi (v),  \\
\psi ([e_0 , v]) & = \psi ( \varphi _1 (e_0) v),
\end{align*}
and therefore $\gamma\varphi_2(e_0)\psi(v)=\psi(\varphi_1(e_0)v)$
for all $v\in\RR^3$, that is,
$\gamma\varphi_2(e_0)=\psi\varphi_1(e_0)\psi^{-1}$.

The converse is straightforward.
\end{proof}

Dozias and Mubarakzyanov gave in \cite{D} and \cite{Mu} a classification of four dimensional solvable Lie algebras. We prove below this result to make this article self contained. The proof uses Proposition~\ref{ideal} together with Lemma~\ref{der}.

\begin{thm}\label{classes}
Let $\ggo$ be a four-dimensional solvable real Lie algebra. Then
$\ggo$ is isomorphic to one and only one of the following Lie
algebras: $\RR^4, \; \; \aff(\RR)\times\aff(\RR), \;\; \RR\times
\hh_3, \;\; \RR \times \mathfrak r_3, \;\; \RR\times \mathfrak
r_{3,\lambda}, \;\; |\lambda|\leq 1, \;\; \RR\times \mathfrak
r'_{3,\lambda}, \;\; \lambda\geq 0, \;\;$ or one of the Lie
algebras with brackets given below in the basis $e_i,\; \;i=0,1,2,3$:

\begin{enumerate}

\item [$\mathfrak n_4$: ] $[e_0, e_1]= e_{2}, \; [e_0, e_2]= e_3$;

\item [$\aff(\CC)$: ] $[e_0,e_2]=e_2, \; \; [e_0,e_3]=e_3, \;\; [e_1, e_2]= e_3,
\; \; [e_1,e_3] =- e_2$;

\item [$\mathfrak r _{4}$: ] $[e_0,e_1]=e_1, \;\; [e_0, e_2]=e_1+e_2, \;\;
[e_0,e_3]=e_2+e_3$;

\item [$\mathfrak r _{4, \lambda}$: ] $[e_0, e_1]=e_1, \;\;
[e_0,e_2]=\lambda e_2,\;\; [e_0, e_3]= e_2+\lambda e_3$;

\item [$\mathfrak r _{4, \mu ,\lambda}$: ] $[e_0,e_1]=e_1, \; \;
[e_0, e_2]= \mu e_2, \;\; [e_0,e_3]= \lambda e_3, \;\; \mu
\lambda\neq 0, \; \; -1 < \mu \leq \lambda\leq 1\;$ or $\; -1=\mu \leq \lambda <0 $ ;

\item [$\mathfrak r _{4, \mu ,\lambda}'$: ] $[e_0,e_1]=\mu e_1, \;\;
[e_0,e_2]=\lambda e_2 -e_3, \;\; [e_0,e_3]=e_2 +\lambda e_3, \;\; \mu > 0$;

 \item [$ \mathfrak d _{4} $: ]  $[e_0,e_1]=  e_1,
\; \;[e_0,e_2]=  -e_2,  \; \; [e_1,e_2]=e_3$;

\item [$ \mathfrak d _{4, \lambda} $: ]  $[e_0,e_1]= \lambda e_1,
\;\; [e_0,e_2]= (1 -\lambda)e_2, \;\; [e_0,e_3]=  e_3, \;\; [e_1,e_2]=e_3,
\;\; \lambda \geq \dfrac 12$;

\item [$ \mathfrak d _{4, \lambda}'$: ] $[e_0, e_1]=\lambda e_1
-e_2, \;\; [e_0, e_2]= e_1+\lambda e_2, \;\; [e_0, e_3]= 2 \lambda
e_3, \;\; [e_1,e_2]=e_3, \;\; \lambda \geq 0 $;

\item [$\mathfrak h _4$: ] $[e_0,e_1]=e_1, \; [e_0,e_2]=e_1+e_2, \;
[e_0,e_3]=2e_3, \;\;[e_1,e_2]=e_3$.

\end{enumerate}

Among these, the unimodular algebras  are: $\RR^4, \; \;  \RR\times
\hh_3, \;\; \RR\times \mathfrak r_{3,-1}, \;\;  \RR\times \mathfrak
r'_{3,0},\;\;\mathfrak n_4, \;\; \mathfrak r_{4,-1/2},\linebreak \mathfrak r_{4,\mu,-1-\mu}$ {\scriptsize $(-1<\mu\leq -1/2)$}, $\;\; \mathfrak r'_{4,\mu,-\mu/2},\,\;\;\mathfrak d _{4},\;\; \mathfrak d_{4,0}'.$
\end{thm}

\begin{proof}
In view of Proposition \ref{ideal} there exists a three
dimensional ideal $\mathfrak v$ of $\ggo$ isomorphic to
$\RR^3,\;\mathfrak e(2),\linebreak\mathfrak e (1,1)$ or $\mathfrak h_3$. We will analyze below the different cases.

\subsection{Case $\mathfrak v= \RR^3$.} We introduce first the following
$3\times 3$ real matrices which will be needed in the next paragraphs:
\begin{eqnarray} \label{fi2}
A_1^{\mu, \lambda} &= &
\begin{pmatrix} 1 & 0 & 0 \\
0 & \mu & 0\\
0& 0& \lambda
\end{pmatrix},  \qquad \qquad
A_2^{\lambda}=\begin{pmatrix}
1 & 0 &0 \\
0 & \lambda & 1\\
0& 0& \lambda
\end{pmatrix},    \\ \label{fi5}
A_3&=&\begin{pmatrix}
1 & 1 &0 \\
0 & 1 & 1\\
0& 0& 1
\end{pmatrix}, \qquad  \qquad
A_4^{\mu, \lambda}=\begin{pmatrix}
\mu & 0 &0 \\
0 & \lambda & 1\\
0& -1& \lambda
\end{pmatrix}.
\end{eqnarray}

By assumption, $\ggo=\RR e_0\ltimes_{\varphi}\RR^3$ where
$\varphi(e_0)=\ad(e_0)$. Suppose first that $\varphi(e_0)$ has
real eigenvalues. We have the following possibilities for
$\varphi(e_0)$, where the eigenvalues are ordered such that
$|\lambda _1 |\leq |\lambda _2 |\leq |\lambda _3 |$:
\[ \text{i) }  \varphi(e_0)=\begin{pmatrix}
\lambda_1 & 0 & 0 \\
0 & \lambda_2 & 0\\
0& 0& \lambda_3
\end{pmatrix}, \quad
\text{ii) } \varphi(e_0)=\begin{pmatrix}
\lambda_1 & 0 &0 \\
0 & \lambda_2 & 1\\
0& 0& \lambda_2
\end{pmatrix}, \quad
\text{iii) } \varphi(e_0)=
\begin{pmatrix}
\lambda & 1 & 0 \\
0 & \lambda & 1\\
0& 0& \lambda
\end{pmatrix}. \]

\bigskip

Case i) $\begin{cases}
\lambda_i=0,\, i=1,2,3, \text{ then } \ggo\cong\RR^4;\\
\lambda_1=0,\;\lambda_3\neq 0,\text{ then } \ggo \cong \mathfrak r
_{3,\lambda}\times\RR; \text{ where } \lambda=
\frac{\lambda_2}{\lambda_3};\\
\lambda_1\lambda _2\lambda_3 \neq 0 \text{ then } \ggo \cong
\RR\ltimes_{\varphi_1} \RR^3,\; \;
\varphi_1(e_0)=A_1^{\mu,\lambda} \text{ as shown in \eqref{fi2},
that is, } \ggo\cong\mathfrak r_{4,\mu ,\lambda}.
\end{cases} $

\medskip

The last isomorphism in Case i) follows by dividing $e_0$ by
$\lambda _3$ and by reordering suitably the basis
$\{e_1,e_2,e_3\}$ of $\RR^3$, we may assume that
$-1\leq\mu\leq\lambda\leq 1$.

\bigskip

Case ii) $\begin{cases}
\lambda_1= \lambda _2=0 \text{ then } \ggo \cong  \RR \times\mathfrak h_3 ;\\
\lambda_1=0, \; \lambda_2\neq 0 \text{  then } \ggo \cong
\RR\times \mathfrak r_3;\\
\lambda _1 \neq 0,  \text{ then } \ggo \cong \RR \ltimes _{\varphi
_{2}^{\lambda}} \RR^3,
\;\;\varphi_2^{\lambda}(e_0)=A_2^{\lambda}\text{ as shown in }
\eqref{fi2}, \text{ that is, } \ggo \cong \mathfrak r _{4,
\lambda}.
\end{cases}$

\bigskip

Case iii) $\begin{cases}
\lambda = 0, \text{ then } \ggo \cong \mathfrak n_4;\\
\lambda \neq 0, \text{ then } \ggo \cong \RR
\ltimes_{\varphi_{3}}\RR^3, \;\varphi_{3}(e_0)=A_3\text{ as shown
in } \eqref{fi5}, \text{ that is, } \ggo \cong \mathfrak r _4.
\end{cases}$

\medskip

The last isomorphism in case iii)  follows by taking $e_0 / \lambda$.

\

In case $\varphi(e_0)$ has only one real eigenvalue, $\mu$, then
we may assume that $\varphi(e_0)=A_4^{\mu, \lambda}$ as in
\eqref{fi5} and we have:
\[   \hspace{-7.5cm}\begin{cases}
 \mu=0, \text{ then } \ggo \cong
\RR \times \mathfrak r'_{3,\lambda};  \\  \vspace{-.3cm} \\
\mu\neq 0, \text{ then } \ggo \cong \mathfrak r '_{4, \mu ,
\lambda}, \; \mu > 0.\end{cases}   \] Observe that the last
isomorphism follows by changing $e_0$ by $-e_0$.

\subsection{Case $\mathfrak v= \mathfrak e(2)$.}

Assume that $\ggo=\RR e_0\ltimes_{\varphi}\mathfrak e(2)$ where $\varphi (e_0)=\ad(e_0)\in
\Der\mathfrak e(2)$ is as in \eqref{der1}. Then setting $e_0'=e_0-be_1+de_2-ce_3\;$, it follows that
\[ [e_0', e_1]=0, \qquad [e_0', e_2]=a e_2,  \quad [e_0', e_3]=ae_3; \]
therefore, $\ggo \cong \RR \times \mathfrak e (2)= \RR \times
\mathfrak r '_{3,0}$ or $\ggo\cong \aff(\CC)$ depending on $a=0$
or $a\neq 0$, respectively.

\subsection{Case $\mathfrak v= \mathfrak e(1,1)$.}

Assume that $\ggo=\RR e_0 \ltimes_{\varphi} \mathfrak  e(1,1)$
where $\varphi (e_0)=\ad(e_0)\in \Der\mathfrak  e(1,1)$ is as in \eqref{der2}.
Let $e_0'=e_0-ae_1+c e_2 -de_3$, then
\[ [e_0', e_1]=0, \quad [e_0', e_2]=0, \qquad [e_0', e_3]=(a+b) e_3;\]
therefore, $\ggo \cong \RR \times \mathfrak e (1,1)= \RR \times
\mathfrak r _{3,-1}$ or $\ggo\cong\aff(\RR)\times\aff(\RR)$
depending on $a+b=0$ or $a+b\neq 0$, respectively.

\subsection{Case $\mathfrak v= \mathfrak h_3 $. }
Assume that $\ggo\cong \RR e_0\ltimes_{\varphi}\mathfrak h_3$
where $\varphi(e_0)=\ad(e_0)$ is given by
\[ \begin{pmatrix}
A & \begin{matrix} 0 \\  0 \end{matrix}  \\
\, b \;\; \; c & \text{tr} A
\end{pmatrix} \] (see \eqref{der3}).
We may assume that $b=c=0$. In fact, setting $e_0'=e_0-c e_1+b
e_2$ it turns out that $\ad(e_0')$ is given by
\begin{equation}
 \begin{pmatrix}
 A & \begin{matrix} 0 \\  0 \end{matrix}  \\
\, 0 \;\; \; 0 & \text{tr} A
\end{pmatrix}.
\end{equation}

Assume first that $A$ has two real eigenvalues $\gamma,\;\beta$;
then $A$ takes the form
\[ \text{i) } A=\begin{pmatrix} \gamma & 0\\
                0& \beta
\end{pmatrix}, \hspace{1.5cm} \text{or} \hspace{1.5cm}
\text{ii) } A=\begin{pmatrix} \gamma & 1\\
               0& \gamma
\end{pmatrix}. \]
Observe that in all cases, once we change the basis $e_1, e_2$ to
$e_1',e_2'$, we must set $e_3'=[e_1',e_2']$ in order to obtain a
Lie algebra isomorphism.

\begin{equation*} \hspace{-1cm}
\text{ Case i)} {\begin{cases}\gamma=\beta=0,\text{ then }
\ggo\cong \RR \times \mathfrak h _3;\\
\gamma=-\beta \neq 0,\text{ then }
\ggo \cong \mathfrak d _{4};\\
\gamma+\beta \neq 0, \text{ then } \ggo \cong \mathfrak d _{4,
\lambda}, \;\;  \lambda = \dfrac {\gamma}{\gamma +\beta};
\end{cases} }  \qquad
\text{Case ii)} {\begin{cases} \gamma=0, \text{ then } \ggo \cong \mathfrak n_4;\\
\gamma \neq 0, \text{ then } \ggo \cong \mathfrak h_4.
\end{cases}}
\end{equation*}

We show that $\mathfrak d_{4,\lambda} \cong \mathfrak
d_{4,1-\lambda}$. This follows by changing the basis $e_i, \; 0\leq
i\leq 3,$ to the basis $e_i', \; 0\leq i\leq 3,$ where:
\[ e_0'=e_0, \quad e_1'=e_2, \quad e_2' =e_1, \quad e_3'= -e_3. \]
Therefore, we may assume that $\lambda \geq 1/2$.

In case ii), $\gamma \neq 0$, in order to show that $\ggo
\cong\mathfrak h_4$ one has to start with
$e_0'=\dfrac1{\gamma}e_0$, then take $e_1',e_2'\in
\text{span}\{e_1, e_2 \}$ such that
\[ \ad(e_0')= \begin{pmatrix} 1 & 1\\
                0& 1
\end{pmatrix} \]
with respect to $\{e_1',e_2'\}$ and $e_3'=[e_1',e_2']$.

If $A$ has no real eigenvalues, then $\ad(e_0')$ takes the form
\[ \begin{pmatrix} \lambda & 1& 0\\
                   -1 & \lambda & 0\\
                    0 & 0 & 2\lambda
\end{pmatrix}, \]
and we conclude that $\ggo\cong\mathfrak d_{4,\lambda}'$. Hence, we have shown so far that any four dimensional solvable Lie algebra is isomorphic to one of those listed in the statement of
the theorem. It remains to show that they are pairwise non isomorphic.

\medskip

\subsection{Isomorphism classes}

In Table \ref{comm}, we list the four dimensional solvable Lie algebras according to their commutator. After that, we proceed to distinguish them up to isomorphism.

\begin{table}
\begin{center}
\begin{tabular}{|l|l|}\hline
${\begin{array}{c} \vspace{-.35cm}\\ \;\;[\ggo,\ggo]\\
\vspace{-.35cm}
\end{array}}$ & $\;\;\;\;\ggo $    \\
\hline \hline ${\begin{array}{c} \vspace{-.35cm}\\ \{ 0 \} \\
\vspace{-.35cm}
\end{array} }$ & $\RR^4$\\
\hline
${\begin{array}{c} \vspace{-.35cm}\\ \RR
\\ \vspace{-.35cm}
\end{array} }$ &
$\RR \times \mathfrak h_3 \, , \;\; \RR \times \mathfrak r _{3,0}$ \\
\hline
${\begin{array}{c} \vspace{-.35cm}\\ \RR^2, \;\mathfrak z=\{ 0\}\\
\vspace{-.35cm}
\end{array} }$ &
$\aff(\RR)\times\aff(\RR), \;\; \aff(\CC),\;\; \mathfrak d _{4,1}$\\
\hline
${\begin{array}{c} \vspace{-.35cm}\\ \RR^2, \;\mathfrak z\neq \{ 0\} \\
 \vspace{-.35cm}
\end{array} }$ &
$\RR\times\mathfrak r_3, \;\; \RR\times\mathfrak r_{3, \lambda}$
({\scriptsize $|\lambda| \leq 1, \;\lambda \neq 0$}),
$\;\; \RR\times \mathfrak r'_{3,\lambda}$ ({\scriptsize $\lambda \geq 0$}),
$\;\; \mathfrak r_{4,0}, \;\; \mathfrak n _4$\\

\hline
${\begin{array}{c} \vspace{-.35cm}\\ \RR^3 \\
\vspace{-.35cm}
\end{array} }$ &
$\mathfrak r_{4}\, ,\;\; \mathfrak r_{4, \lambda}$ ({\scriptsize $\lambda\neq 0$})  ,\; $\mathfrak r_{4,\mu,\lambda}\;$
({\scriptsize $\mu \lambda \neq 0, \; -1 \leq \mu \leq \lambda \leq 1)$},$\;\;
\mathfrak r'_{4,\mu , \lambda}$ ({\scriptsize $\mu >0$})  \\
\hline
${\begin{array}{c} \vspace{-.35cm}\\ \mathfrak h_3 \\
\vspace{-.35cm}
\end{array} }$ &
$\mathfrak d_4, \;\; \mathfrak d_{4, \lambda}$ ({\scriptsize $\lambda \neq 1,\;
\lambda \geq 1/2 $}), $\;\; \mathfrak d'_{4, \lambda}$ ({\scriptsize $\lambda \geq 0$}),
$\;\; \mathfrak h_4$  \\
\hline
\end{tabular}
\bigskip
\caption{}\label{comm}
\end{center}
\end{table}

\smallskip

$\bullet\, [\ggo,\ggo ]= \RR$: $\RR \times \mathfrak h_3$ is
nilpotent but $\RR \times \mathfrak r_{3,0}$ is not, therefore
they are not isomorphic.

\medskip

$\bullet\,[\ggo,\ggo]=\RR^2,\;\;\mathfrak z =\{ 0\}$: Both
$\aff(\RR)\times\aff(\RR)$ and $\mathfrak d_{4,1}$ are completely
solvable\footnote{Recall that a solvable Lie algebra $\ggo$ is {\em completely solvable} when $\ad(x)$ has real eigenvalues for all $x\in\ggo$.} and therefore not isomorphic to $\aff(\CC),$ which is not
completely solvable. The unimodular kernel of $\aff(\RR) \times
\aff(\RR)$ (resp. $\mathfrak d_{4,1}$) is $\mathfrak r_{3, -1}$
(resp. $\mathfrak h_3$), hence $\aff(\RR)\times\aff(\RR)$ is not
isomorphic to $\mathfrak d _{4,1}$.

\medskip

$\bullet\,[\ggo,\ggo]= \RR^2, \;\; \mathfrak z \neq \{ 0\}$: If
$\ggo=\RR\times\mathfrak r_3, \;\; \RR \times \mathfrak
r_{3,\lambda}\, (|\lambda|\leq 1,\,\lambda \neq 0)$ or
$\RR\times\mathfrak r'_{3,\lambda}\, (\lambda\geq 0)$ then
$\mathfrak z \cap[\ggo,\ggo]=\{ 0\}$, while $\mathfrak z \cap
[\ggo,\ggo] \neq \{ 0\}\;$ when $\; \ggo = \mathfrak r_{4,0}$ or
$\mathfrak n _4$. Also $\ggo = \RR \times \mathfrak r_3, \;\; \RR
\times \mathfrak r_{3, \lambda}\, (|\lambda|\leq 1,\,\lambda \neq
0)$ and $\RR \times \mathfrak r'_{3, \lambda}\,(\lambda\geq 0)$
are not pairwise isomorphic since $\mathfrak r_3, \; \mathfrak
r_{3, \lambda}$ and $\mathfrak r'_{3,\lambda}$ are 3-dimensional
non isomorphic Lie algebras. On the other hand, $\mathfrak n_4$ is
nilpotent but $\mathfrak r_{4,0}$ is not, hence they are not
isomorphic.

\medskip

$\bullet\,[\ggo,\ggo]=\RR^3$: $\mathfrak r _{4},\; \mathfrak r
_{4, \lambda}\, (\lambda \neq 0), \; \mathfrak r_{4,\mu, \lambda},
\; \mathfrak r'_{4,\mu,\lambda}$. In this case, it follows from
Lemma \ref{R3} that any pair of Lie algebras belonging to
different families can not be isomorphic. The last family consists
of non completely solvable Lie algebras.

The fact that two Lie algebras $\mathfrak r_{4,\lambda},\; \lambda
\neq 0$, and $\mathfrak r _{4,\lambda'},\;\lambda'\neq 0$, are
isomorphic if and only if $\lambda =\lambda'$ follows by applying
Lemma \ref{R3}.

Let us show that if $\mathfrak r_{4,\mu,\lambda},\;-1< \mu \leq
\lambda \leq 1, \; \mu \lambda \neq 0$, is isomorphic to
$\mathfrak r_{4,\mu',\lambda'},\;-1 < \mu'\leq \lambda' \leq
1,\; \mu'\lambda'\neq 0$, then $\mu = \mu '$ and $\lambda =
\lambda '$. From Lemma \ref{R3}, there exists $\gamma \neq 0$ such
that the sets of eigenvalues $\{ 1, \mu , \lambda \}$ and $\{
\gamma, \gamma \mu ', \gamma \lambda '\}$ must coincide. If
$\gamma =1$ the desired assertion follows from  $\mu \leq \lambda$
and $\mu '\leq \lambda '$. If $\gamma=\mu$ then either $\gamma \mu'=1$ or $\gamma \lambda '=1$, hence $\mu '=1$ or $\lambda '=1$, therefore $\gamma=1$ and again
this implies $\mu=\mu ', \; \; \lambda =\lambda '$. The case $\gamma =\lambda$ is proved in a similar way.

Let us show that if $\mathfrak r_{4,-1,\lambda},\;-1 \leq \lambda <0 $, is isomorphic to
$\mathfrak r_{4,-1,\lambda'},\;-1 \leq \lambda' <0$, then  $\lambda =
\lambda '$. We apply
Lemma~\ref{R3} again to obtain that there exists $\gamma\neq 0$
such that $\{ 1, -1 , \lambda \}$ and $\{
\gamma, -\gamma , \gamma \lambda '\}$ must coincide. We cannot have $\gamma=-1$, since this would imply
$\lambda =-\lambda '$, a contradiction, since both, $\lambda$ and $\lambda'$
are negative. If $\gamma=\lambda$, then $-\gamma= 1$ and $-1=\gamma \lambda '
=\lambda \lambda '>0$, a contradiction. Thus $\gamma =1$ and
 $\lambda =\lambda '$.

If $\mathfrak r_{4,\mu,\lambda},\;-1< \mu \leq
\lambda \leq 1, \; \mu \lambda \neq 0$, were isomorphic to
$\mathfrak r_{4,-1,\lambda '},\;-1 \leq \lambda ' <0 $, then
Lemma~\ref{R3} would imply that that there exists $\gamma\neq 0$
such that $\{ 1, \mu , \lambda \}=\{
\gamma, -\gamma , \gamma \lambda '\}$. If $\gamma=1$ then $\mu=-1$, which is impossible. On the other hand, $\gamma = \mu$ implies $-\gamma =1$ or $\gamma \lambda '=1$, hence $\mu=\gamma =-1$, a contradiction. The case $\gamma = \lambda $ is similar; therefore, the above Lie
algebras are not isomorphic.

Assume now that $\mathfrak r'_{4,\mu,\lambda},\;\mu >0$, is
isomorphic to $\mathfrak r'_{4,\mu',\lambda'}, \;\mu'>0$, we must
show that $\mu =\mu'$ and $\lambda = \lambda'$. We apply
Lemma~\ref{R3} again to obtain that there exists $\gamma\neq 0$
such that $\mu =\gamma \mu'$ and $\lambda \pm i = \gamma
(\lambda'\pm i)$. It follows from the second equality that $\gamma
= \pm 1$, and the first equality implies $\gamma =1$, since both
$\mu$ and $\mu'$ are positive. Therefore, $\mu =\mu'$ and $\lambda
= \lambda'$, as claimed.

\medskip

$\bullet\,[\ggo ,\ggo ]= \mathfrak h_3$: The Lie algebras
$\mathfrak d_4, \; \mathfrak d_{4,\lambda}\; (\lambda \geq 1/2, \,
\lambda \neq 1)$, and $\mathfrak d'_{4, \lambda}, \; \mathfrak h_4$
are distinguished by $\ggo / \mathfrak z( [\ggo,\ggo])$, as the
following table shows:

\medskip

\begin{center}
\begin{tabular}{|c||c|c|c|c|} \hline
$ \;\;\ggo$ & $\mathfrak d_4$ & $\mathfrak d_{4, \lambda}, \;$
\scriptsize{$\begin{cases} \lambda \geq 1/2  \\\lambda \neq 1
\end{cases}$} &
$\mathfrak d'_{4, \lambda}, \;$ \scriptsize{$\lambda \geq 0$} & $\mathfrak h_4$\\
\hline
${\begin{array}{c} \vspace{-.2cm}\\ \;\;\ggo/\mathfrak z( [\ggo , \ggo])\\
\vspace{-.2cm}
\end{array}}$ & $\mathfrak r_{3, -1}$ &
$ \mathfrak r_{3,-1 + 1/\lambda }$ & $
\mathfrak r'_{3,\lambda }$ & $ \mathfrak r_{3}$  \\
\hline
\end{tabular}
\end{center}
\end{proof}

\begin{remarks} (i) In \cite{D-S} it was proved that $\mathfrak d'_{4, \lambda}, \;
\lambda \geq 0$, are all non-isomorphic. Observe that $\ggo
_{\lambda}$ in \cite{D-S} corresponds to $\mathfrak d'_{4, 1
/\lambda}$ for $\lambda \neq 0$ (resp. $\mathfrak d_{4, 1/2}$ for
$\lambda =0$).

\noindent (ii) We observe that $\aff(\CC)$ is the Lie algebra of
the group of affine motions of the complex line, which is
isomorphic to the complexification of $\aff(\RR)$ looked upon as a
real Lie algebra. Also, $\mathfrak r_{4,1,1}$ is the Lie algebra of a
solvable Lie group which acts simply and transitively on the real
hyperbolic space $\RR H^4$ and $\mathfrak d_{4,1/2}$ is the Lie
algebra of a solvable Lie group which acts simply and transitively
on the complex hyperbolic space $\CC H^2.$
\end{remarks}

\

\section{Product structures on four dimensional solvable Lie algebras }

\subsection{Basic definitions}
An {\em almost product structure} on a Lie algebra $\ggo$ is a
linear endomorphism $E:\ggo \longrightarrow \ggo$ satisfying $E^2=\,$Id (and
not equal to $\pm$Id). It is said to be {\em integrable} if \begin{equation}
\label{integrable2} E[x,y]=[Ex,y]+[x,Ey]-E[Ex,Ey] \quad \text{for all } x,y\in
\ggo.  \end{equation}
An integrable almost product structure will be called a
{\em product structure}.

An almost product structure on $\ggo$ gives rise to a decomposition of $\ggo$ into
\begin{equation}\label{decomposition}
\ggo= \ggo_{+} \oplus \ggo_{-},\;\;\; E|\ggo_{+}=\text{Id},\; E|\ggo_{-}=-\text{Id}.
\end{equation}
The integrability of $E$ is equivalent to $\ggo_{+}$ and $\ggo_{-}$ being subalgebras.
When $\dim\ggo_+=\dim\ggo_-$, the product structure $E$ is called a {\em paracomplex structure}.

Three Lie algebras $(\ggo,\ggo_+,\ggo_-)$ form a {\em double Lie algebra} if
$\ggo_+$ and $\ggo_-$ are Lie subalgebras of $\ggo$ and
$\ggo=\ggo_+\oplus\ggo_-$ as vector spaces. This will be denoted by $\ggo=\ggo_+ \bowtie \ggo_-$.
Observe that a double Lie algebra $(\ggo,\ggo_+,\ggo_-)$ gives a product
structure $E:\ggo\longrightarrow\ggo$ on $\ggo$, where $E|{\ggo_+}=\,$Id and
$E|{\ggo_-}=- $Id. Conversely, a product structure on the Lie algebra $\ggo$
gives rise to a double Lie algebra $(\ggo,\ggo_+,\ggo_-)$, where $\ggo_{\pm}$
is the eigenspace associated to the eigenvalue $\pm 1$ of $E$.
The notion of double Lie algebra is
a natural generalization of that of semidirect product. We will denote $\ggo=\ggo_+\ltimes\ggo_-$ the semidirect product of $\ggo_+$ and $\ggo_-$ where $\ggo_-$ is an ideal of $\ggo$, that is, there is a split exact sequence
\[ 0\longrightarrow \ggo_-\longrightarrow \ggo \longrightarrow \ggo_+ \longrightarrow 0.\]

Product structures or, equivalently, double Lie algebras, were used in several contexts (see \cite{AS}, \cite{LW}). Important examples of double Lie algebras are Manin triples and complex product structures.

\subsection{Paracomplex structures }
It is the main goal of this subsection to determine all 4-dimensional solvable Lie algebras admitting  paracomplex structures. We will give realizations of the Lie algebras obtained in Theorem~\ref{classes} as double Lie algebras with subalgebras of dimension 2 when such a structure exists (see Table \ref{pc}), or prove the non existence otherwise. It turns out that among all four dimensional solvable Lie algebras there is only one family, whose commutator ideal is $\hh_3$, not admitting any paracomplex structure (Theorem~\ref{paracomplex}). Since there are only two non-isomorphic two-dimensional Lie algebras: $\RR^2$ and $\aff(\RR)$, the possible decompositions $\ggo_+ \bowtie \ggo_-$ are $\RR^2 \bowtie \RR^2,\; \RR^2 \bowtie \aff(\RR)$ and $\aff(\RR) \bowtie \aff(\RR)$.

\begin{table}
\begin{center}
\small{
\begin{tabular}{|c|c|c|c|}    \hline
$\begin{array}{c} \vspace{-.35cm}\\ \ggo
\\ \vspace{-.35cm} \end{array}$ & $\RR^2 \bowtie \RR^2$ & $\aff(\RR) \bowtie \RR^2$ & $\aff(\RR) \bowtie \aff(\RR) $\\ \hline
 $\begin{array}{c} \vspace{-.35cm}\\ \RR^4 \\ \vspace{-.35cm} \end{array}$ & $\langle e_0,e_1\rangle\times \langle e_2,e_3\rangle$ & no & no \\ \hline

$\begin{array}{c} \vspace{-.35cm}\\ \aff(\RR)\times\aff(\RR)
\\ \vspace{-.35cm} \end{array}$ & $\langle e_0,e_1\rangle\ltimes \langle e_2,e_3\rangle$ & $\langle e_1+e_3,e_2\rangle\bowtie\langle e_0,e_1\rangle $ & $\langle e_0,e_3\rangle\times \langle e_1,e_2\rangle$\\ \hline

$\begin{array}{c} \vspace{-.35cm}\\ \RR\times\hh_3
\\ \vspace{-.35cm} \end{array}$ &$\langle e_0,e_2\rangle\ltimes \langle e_1,e_3\rangle$ & no & no \\ \hline

$\begin{array}{c} \vspace{-.35cm}\\ \RR \times \mathfrak r_3
\\ \vspace{-.35cm} \end{array}$  & $\langle e_0,e_1\rangle\ltimes \langle e_2,e_3\rangle$ & $\langle e_1,e_2\rangle\bowtie\langle e_0,e_3\rangle$ & no \\ \hline

$\begin{array}{c} \vspace{-.35cm}\\ \RR\times \mathfrak r_{3,\lambda},\,\lambda\neq 0 \\ \vspace{-.35cm} \end{array}$ & $\langle e_0,e_1\rangle\ltimes \langle e_2,e_3\rangle$ & $\langle e_1,e_2\rangle\ltimes\langle e_0,e_3\rangle$ & $\langle e_0+e_1,e_2\rangle\bowtie \langle e_1-\lambda e_0,e_3\rangle$ \\ \hline

$\begin{array}{c} \vspace{-.35cm}\\ \RR\times \mathfrak r_{3,0}
 \\ \vspace{-.35cm} \end{array}$ & $\langle e_0,e_1\rangle\ltimes \langle e_2,e_3\rangle$ & $\langle e_1,e_2\rangle\times\langle e_0,e_3\rangle$ & no \\ \hline

 $\begin{array}{c} \vspace{-.35cm}\\ \RR\times \mathfrak r'_{3,\lambda}
 \\ \vspace{-.35cm} \end{array}$ & $\langle e_0,e_1\rangle\ltimes \langle e_2,e_3\rangle$& no & no \\ \hline

$\begin{array}{c} \vspace{-.35cm}\\ \mathfrak n_4
\\ \vspace{-.35cm} \end{array}$ &$\langle e_0,e_3\rangle\bowtie \langle e_1,e_2\rangle$ &no& no \\ \hline

$\begin{array}{c} \vspace{-.35cm}\\ \aff(\CC)
 \\ \vspace{-.35cm} \end{array}$ & $\langle e_0,e_1\rangle\ltimes \langle e_2,e_3\rangle$  & $\langle e_0, e_2\rangle\bowtie \langle e_0- e_3, e_1 + e_2\rangle$ & no\\ \hline

$\begin{array}{c} \vspace{-.35cm}\\ \mathfrak r _{4}
 \\ \vspace{-.35cm} \end{array}$ &no& $\langle e_0, e_1\rangle\bowtie \langle e_2, e_3\rangle$ & no \\ \hline

$\begin{array}{c} \vspace{-.35cm}\\ \mathfrak r _{4, \lambda},\; \lambda\neq 0
\\ \vspace{-.35cm} \end{array}$ & no & $\langle e_0, e_1\rangle\ltimes \langle e_2, e_3\rangle $ & $\langle e_0,e_1\rangle\bowtie\langle e_0+{\lambda}e_3,e_2\rangle$\\ \hline

$\begin{array}{c} \vspace{-.35cm}\\ \mathfrak r _{4, 0}
 \\ \vspace{-.35cm} \end{array}$ & $\langle e_0, e_2\rangle\bowtie\langle e_1, e_3\rangle$ & $\langle e_0, e_1\rangle\ltimes \langle e_2, e_3\rangle$ & no \\ \hline

$\begin{array}{c} \vspace{-.35cm}\\ \mathfrak r _{4, \mu ,\lambda}
 \\ \vspace{-.35cm} \end{array}$ &no& $\langle e_0, e_1\rangle\ltimes \langle e_2, e_3\rangle $&$\langle e_0-e_1, e_2\rangle\bowtie \langle e_0+e_1, e_3\rangle$ \\ \hline

$\begin{array}{c} \vspace{-.35cm}\\ \mathfrak r _{4, \mu ,\lambda}'
 \\ \vspace{-.35cm} \end{array}$ & no & $\langle e_0, e_1\rangle\ltimes \langle e_2, e_3\rangle$ & no \\ \hline

 $\begin{array}{c} \vspace{-.35cm}\\ \mathfrak d _{4}
\\ \vspace{-.35cm} \end{array}$ & no & $\langle e_0, e_1\rangle\ltimes \langle e_2, e_3\rangle$ & $\langle e_0+e_2, e_1-e_3\rangle  \bowtie \langle e_0-e_2, e_1+e_3\rangle$
\\ \hline

 $\begin{array}{c} \vspace{-.35cm}\\ \mathfrak d _{4, \lambda}, \lambda \neq 1 \\ \vspace{-.35cm} \end{array} $ & no& $\langle e_0, e_1\rangle\ltimes \langle e_2, e_3\rangle$&
$\langle e_0, e_3\rangle \bowtie \langle e_0 +{\lambda} e_2, {(1-\lambda)} e_1 +{\lambda}e_3\rangle$ \\ \hline

$\begin{array}{c} \vspace{-.35cm}\\ \mathfrak d _{4, 1}
\\ \vspace{-.35cm} \end{array}$ & $\langle e_0, e_2\rangle\ltimes \langle e_1, e_3\rangle$ & $\langle e_0, e_1\rangle\ltimes \langle e_2, e_3\rangle$ & $\langle e_0,e_1\rangle\bowtie\langle e_0+e_2,e_3\rangle$
 \\ \hline

$\begin{array}{c} \vspace{-.35cm}\\  \mathfrak d _{4, \lambda}'
\\ \vspace{-.35cm} \end{array}$&no&no& no\\ \hline

 $\begin{array}{c} \vspace{-.35cm}\\ \mathfrak h _4
 \\ \vspace{-.35cm} \end{array}$ & no & $\langle e_0, e_1\rangle\bowtie \langle e_2, e_3\rangle$ &
$\langle e_0, e_3\rangle\bowtie \langle e_0-e_2, e_1-e_3\rangle$ \\ \hline

\end{tabular}}
\end{center} \vspace{.2cm}
\caption{Paracomplex structures on four dimensional solvable Lie algebras}\label{pc}
\end{table}

By simple computations one can verify that the decompositions given in Table~\ref{pc} satisfy the required properties. We prove below the non existence results.

\begin{prop} Let $\ggo$ be a Lie algebra with an abelian commutator ideal $\ggo'$ of codimension $1$. Then any abelian subalgebra of dimension $n> 1$ is contained in $\ggo'$.
\end{prop}

\begin{proof} In this case there is $e_0\in\ggo$ such that $\ad (e_0)$ is an isomorphism of $\ggo'$. Let $\hh$ be an abelian subalgebra of $\ggo, \;
\dim \hh >1$, and let $x,y\in\hh$ linearly independent. If $x=a_0e_0+x',\,y=b_0e_0+y'$ with $a_0,b_0\in\RR$ and $x',y'\in\ggo'$, then
\[ 0=[x,y]=[e_0,a_0y'-b_0x']. \]
This implies that $a_0y'-b_0x'=0$, that is, $a_0y-b_0x=0$ and hence $a_0=b_0=0$. Therefore, $x,y\in\ggo'$, as asserted.
\end{proof}

The previous result together with Table~\ref{comm} imply
\begin{cor} \label{no21}
The Lie algebras $\mathfrak r_{4} ,\; \mathfrak r_{4, \lambda}$ {\scriptsize ($\lambda\neq 0$)}, $\mathfrak r_{4,\mu,\lambda},\;\mathfrak r'_{4,\mu , \lambda}$ do not admit a decomposition of type $\RR^2\bowtie\RR^2$.
\end{cor}

\smallskip

\begin{prop}[\cite{P}] \label{nilpotent} If $\ggo$ is a Lie algebra which admits a decomposition $\ggo = \ggo_+ \bowtie \ggo_-$ with $\ggo_+$ and $\ggo_-$ abelian subalgebras, then $\ggo$ is 2-step solvable (i.e., $\ggo'$ is abelian).
\end{prop}

\begin{proof} If $\ggo = \ggo_+ \bowtie \ggo_-$ with $\ggo_+$ and $\ggo_-$ abelian then  $[(x_1,x_2), (y_1,y_2)]$ is determined by $[(x_1,0), (0,y_2)]= (\alpha(x_1,y_2), \beta(x_1,y_2))$ where $\alpha$ and $\beta$ denote the components on $\ggo_+$ and $\ggo_-$ respectively.  Since the bracket on $\ggo$ satisfies the Jacobi identity one obtains
\begin{enumerate}
\item $\alpha(x_1,\beta(y_1,z_2))= \alpha(y_1,\beta(x_1,z_2)),$
\item $\beta(\alpha(z_1,y_2),x_2)= \beta(\alpha(z_1,x_2),y_2),$
\item $\beta(x_1,\beta(y_1,z_2))= \beta(y_1,\beta(x_1,z_2)),$
\item $\alpha(\alpha(z_1,y_2),x_2))=\alpha(\alpha(z_1,x_2),y_2))$.
\end{enumerate}
Now, using the above relations one can show that
\[\alpha(\alpha(x_1,y_2),\beta(u_1,v_2))=\alpha(\alpha(u_1,v_2),\beta(x_1,
y_2))\]
and
\[\beta(\alpha(x_1,y_2),\beta(u_1,v_2))=\beta(\alpha(u_1,v_2),\beta(x_1,
y_2)).\]
But the above relations immediately imply
\[ [[(x_1,0),(0,y_2)],[(u_1,0),(0,v_2)]]=0 \]
and the assertion follows.
\end{proof}

The above proposition together with Table~\ref{comm} imply

\begin{cor}\label{no22}
The Lie algebras $\dd_{4},\; \dd_{4, \lambda}$ {\scriptsize ($\lambda\neq 1$)}, $\dd'_{4,\lambda},\;\hh_{4}$ do not admit a decomposition of type $\RR^2\bowtie\RR^2$.
\end{cor}

\smallskip

\begin{lem}\label{sub-aff}
The Lie algebras $\RR^4,\, \RR\times\hh_3,\,\nn_4$ and $\RR\times\rr'_{3,\lambda}$ do not contain $\aff(\RR)$ as a subalgebra. Hence, these Lie algebras do not admit decompositions of type $\aff(\RR)\bowtie\RR^2$ or $\aff(\RR)\bowtie\aff(\RR)$.
\end{lem}

\begin{proof}
Since $\RR^4,\; \RR\times\hh_3$ and $\mathfrak n_4$ are nilpotent, they cannot have subalgebras isomorphic to $\aff(\RR)$. Let us show next that the same holds for $\ggo:=\RR\times\rr'_{3,\lambda}$. In fact, assume that there exist $x,y \in \ggo$ such that $[x,y]=y$. Then $y \in \ggo'$. If $x\in \langle e_0,e_2,e_3\rangle$ then $y=0$, thus assume that $x=e_1 +u$. So $[x,y]=[e_1, y]=y$ implies that $y=0$ since $\ad(e_1)$ has no real eigenvalues in $\ggo'$.
\end{proof}

\smallskip

\begin{prop}\label{nose}
The Lie algebras $\RR \times \rr_3, \RR\times\rr_{3,0},\, \aff(\CC),\, \rr_4,\, \rr_{4,0}$ and $\rr'_{4,\mu,\lambda}$ do not admit a decomposition of type $\aff(\RR)\bowtie\aff(\RR)$.
\end{prop}
\begin{proof}
 Let $\ggo:= \RR \times \rr_3$ and $\hh$  a subalgebra of $\ggo$ isomorphic to $\aff(\RR)$. Then $\hh$ has a basis of the form $\{e_1 + u, e_2\}$ with $u \in \langle e_0,e_3\rangle$. Thus, any decomposition of $\RR \times \rr_3$ of the form $\hh_1 +\hh_2$ with $\hh_1 \simeq \aff(\RR) \simeq \hh_2$ is not direct since $e_2 \in \hh_1 \cap \hh_2$.
If $\ggo=\RR\times\rr_{3,0}$, then $\dim \ggo' =1$ and therefore the assertion
follows.

Assume next that $\ggo \cong \aff(\CC)$. Every subalgebra of $\ggo$ isomorphic to $\aff(\RR)$ is of the form $\la e_0 + u, v\ra $ with $u, v \in \langle e_2,e_3\rangle,$ thus it is contained in the subspace spanned by $\{e_0, e_2, e_3 \}.$ Therefore, $\aff(\CC)$ is not of type $\aff(\RR) \bowtie \aff(\RR)$.

If $\ggo$ is either $\rr_4,\, \rr_{4,0}$ or $\rr'_{4,\mu,\lambda}$, one can show that
$e_1\in\ggo $ belongs to any Lie subalgebra of $\ggo$ isomorphic to $\aff(\RR)$ and thus $\ggo$ cannot be decomposed as $\aff(\RR)\bowtie\aff(\RR)$. We give a proof of this fact in
the case $\ggo= \rr_4$.
Let $\mathfrak u= \la u,v\ra $ be a Lie subalgebra of $\mathfrak r_4$ isomorphic to $\aff(\RR)$, with $[u,v]=v$. If $u=\sum_{i=0}^3a_ie_i,\; v=\sum_{i=0}^3b_ie_i$ with $a_i,\,b_i\in\RR,\,i=0,\ldots,3$,  $b_0=0$ since $v \in [\rr_4, \rr_4]$, then we have
\[ \begin{cases}
    a_0b_1+a_0b_2=b_1,\\
    a_0b_2+a_0b_3=b_2,\\
    a_0b_3=b_3
   \end{cases} \quad \text{which implies} \quad
   \begin{cases}
    b_1(a_0-1)=-a_0b_2,\\
    b_2(a_0-1)=-a_0b_3,\\
    b_3(a_0-1)=0. \end{cases} \]
If $a_0-1\neq0$, then $b_3=0$ and thus $b_2(a_0-1)=0$, which implies $b_2=0$. From this, we have $b_1(a_0-1)=0$, and therefore $b_1=0$, i.e. $v=0$, a contradiction. Hence, $a_0=1$ and then $b_2=b_3=0$. Also $b_1\in\RR\setminus\{0\}$ is arbitrary, and we may take $b_1=1$. So,
\[ u=e_0+a_2e_2+a_3e_3,\quad v=e_1,\]
hence, $e_1\in\rr_4$, as asserted. The proofs of the remaining cases are similar.
\end{proof}

\

\begin{thm}\label{paracomplex}
If $\ggo$ is a four dimensional solvable Lie algebra then $\ggo$ does not admit any paracomplex structure if and only if $\ggo$ is isomorphic to  $\mathfrak d'_{4,\lambda}$
for some $\lambda\geq 0$.
\end{thm}

\begin{proof}
We first show that if $\ggo$ is a Lie algebra in the family $\mathfrak d'_{4,\lambda}$ with $\lambda\geq 0$ then $\ggo$ does not admit a paracomplex structure. Let $\mathfrak u$ be a 2-dimensional Lie subalgebra of $\ggo$ with a basis $\{u,v\}$, where $u=\sum_{i=0}^3a_ie_i,\; v=\sum_{i=0}^3b_ie_i$ with $a_i,\,b_i\in\RR,\,i=0,\ldots,3$

\smallskip

\noindent {\bf Case 1:} $[u,v]=0$ and hence $\mathfrak u\cong\RR^2$. In this case we get that
\[ \begin{cases}
    \lambda(a_0b_1-a_1b_0)+a_0b_2-a_2b_0=0,\\
    \lambda(a_0b_2-a_2b_0)-a_0b_1+a_1b_0=0,\\
    2\lambda(a_0b_3-a_3b_0)+a_1b_2-a_2b_1=0.
   \end{cases} \]
From the first two equations we arrive at
$(a_0b_2-a_2b_0)(\lambda^2+1)=0$, and therefore $a_0b_2-a_2b_0=0$,
which in turn implies $a_0b_1-a_1b_2=0$. Summing up, we have
\[(a_0,a_1,a_2)\times(b_0,b_1,b_2)=(-2\lambda(a_0b_3-a_3b_0),0,0)\]
and hence
\[ \begin{cases}
    2\lambda a_0(a_0b_3-a_3b_0)=0,\\
    2\lambda b_0(a_0b_3-a_3b_0)=0.
   \end{cases} \]
We have two cases:

\noindent (i) $\lambda=0$. Then $(a_0,a_1,a_2)\times(b_0,b_1,b_2)=(0,0,0)$ and therefore $(b_0,b_1,b_2)=\beta(a_0,a_1,a_2)$, with $\beta\in\RR$. Since $u$ and $v$ are linearly independent, we must have $b_3-\beta a_3\neq0$. Thus, we obtain that
\[ e_3=\frac{1}{b_3-\beta a_3}(v-\beta u)\in\mathfrak u.\]

\noindent (ii) $\lambda\neq 0$. If we suppose $a_0b_3-a_3b_0\neq 0$, we arrive at a contradiction; thus $a_0b_3-a_3b_0=0$ and $(a_0,a_1,a_2)\times(b_0,b_1,b_2)=(0,0,0)$. As in the previous case, we have that $e_3\in\mathfrak u$.

\medskip

\noindent {\bf Case 2:} $[u,v]=v$ and hence $\mathfrak u\cong\aff(\RR)$. In this case we obtain that $b_0=0$ and
\begin{equation}\label{delta} \begin{cases}
     \lambda a_0b_1+a_0b_2=b_1,\\
     \lambda a_0b_2-a_0b_1=b_2,\\
     2\lambda a_0b_3+a_1b_2-a_2b_1=b_3.
   \end{cases} \end{equation}
Let us observe first that $a_0\neq0$, since otherwise from (\ref{delta}) we obtain that $b_1=b_2=b_3=0$, i.e. $v=0$, a contradiction. Combining now the first two equations from (\ref{delta}), we arrive at \[b_1b_2\left((\lambda a_0-1)^2+a_0^2\right)=0.\] Since clearly $(\lambda a_0-1)^2+a_0^2\neq 0$, we have that $b_1b_2=0$. It is easily seen that this implies $b_1=b_2=0$. Hence, we need only consider now the equation $2\lambda a_0b_3=b_3$, with $b_3\neq0$.

\noindent (i) $\lambda=0$. In this case, we obtain that $b_3=0$, a contradiction. Thus, $\mathfrak d'_{4,0}$ does not have any Lie subalgebra isomorphic to $\aff(\RR)$.

\noindent (ii) $\lambda\neq 0$. Here, since $b_3\neq0$, we have $a_0=\frac{1}{2\lambda}$ and $u$ and $v$ are given by
\[ u=\frac{1}{2\lambda}e_0+a_1e_1+a_2e_2,\quad v=e_3.\]
Note that $e_3\in\mathfrak u$.

\medskip

In all cases, $e_3\in\ggo$ belongs to any 2-dimensional Lie subalgebra of $\ggo$, and hence this Lie algebra cannot be decomposed as a direct sum (as vector spaces) of two 2-dimensional Lie subalgebras.

The theorem follows by observing that the remaining Lie algebras possess paracomplex structures (see Table~\ref{pc}).
\end{proof}

We give next a characterization of the four dimensional solvable Lie algebras which can be decomposed
as a semidirect product of two dimensional subalgebras.

\subsection{Semidirect extensions of $\RR^2$}
Assume that $\ggo$ contains $\RR^2$ as an ideal and that the short exact sequence \[ 0\to \RR ^2 \to \ggo \to \ggo / \RR ^2\to 0 \]
splits. The next result gives a list of the Lie algebras with this property.

\begin{prop} Let $\ggo$ be a four dimensional solvable Lie algebra.

(i) If there is a split exact sequence \[ 0\to \RR ^2 \to \ggo \to \RR^2 \to 0\]
then $\ggo \cong  \RR ^4 , \; \aff(\RR) \times \aff(\RR), \RR \times \mathfrak h _{3}, \;
\RR \times \mathfrak r _{3},
\; \RR \times \mathfrak r _{3, \lambda}, \; \RR \times \mathfrak r ' _{3, \lambda}, \; \aff(\CC)
\, $ or $\, \dd_{4,1}$.

 (ii) If there is a split exact sequence \[ 0\to \RR ^2 \to \ggo \to \aff(\RR) \to 0\]
then $\ggo \cong   \RR \times \mathfrak r _{3, \lambda}, \;  \mathfrak r _{4, \lambda},
\; \mathfrak r _{4, \mu ,\lambda}, \; \mathfrak r' _{4, \mu ,\lambda}, \; \dd _4 \, $ or $\, \dd _{4,\lambda}$.
\end{prop}
\begin{proof}

(i) Table \ref{pc} exhibits decompositions of $\ggo $ as a semidirect product
$\RR^2 \ltimes \RR^2$ in case $\ggo \cong \RR ^4 , \; \aff(\RR) \times
\aff(\RR), \RR \times \mathfrak h _{3}, \; \RR \times \mathfrak r _{3}, \;
\RR \times \mathfrak r _{3, \lambda}, \; \RR \times \mathfrak r ' _{3,
\lambda}, \; \aff(\CC) \, $ or $\, \dd_{4,1}$.
It follows from Corollaries \ref{no21} and \ref{no22} that $\rr _4, \; \rr_{4, \lambda}, \;
 \rr'_{4, \lambda}, \; \lambda \neq 0, \;  \rr_{4, \mu ,\lambda}, \,  \rr'_{4, \mu ,\lambda},
\; \dd_4 , \; \dd_{4, \lambda}, \; \lambda \neq 1, \;  \dd'_{4, \lambda} \,$ and $\, \hh_4\,$
do not admit such a decomposition. It remains to consider the case  $\ggo \cong \nn_4$ or $\rr_{4,0}$.
Assume that $\ggo=\mathfrak a \ltimes \mathfrak b$ with $\mathfrak a \cong \mathfrak b \cong \RR^2$. Then
$\ggo ' \subset \mathfrak b$, hence $\ggo' =\mathfrak b$ since
 in both cases $\ggo '=\RR^2$ (Table~\ref{comm}).

Consider next the case $\ggo\cong \nn _4$, hence $\mathfrak b =\la e_2, e_3 \ra$
and $\aa =\la x, y \ra$ with
$ x= ae_0 +be_1+u , \; y=ce_0 +d e_1 +v , \;ad-bc \neq 0$,
$ u,v \in  \mathfrak b $. We calculate \[
[x,y] = (ad -bc)e_2 + [e_0, av-cu]    \]
which is non zero since the second summand on the right hand side is a multiple of $e_3$. This contradicts
the fact that $\aa \cong \RR ^2$. Therefore, $\nn _4$ does not decompose as $\RR^2
\ltimes \RR^2$.

The case $\ggo \cong \rr_{4, 0}$ is similar. We have $\bb=\la e_1, e_2 \ra $ and
$\aa =\la x, y \ra$ with
$ x= ae_0 +be_3+u , \; y=ce_0 +d e_3 +v , \;ad-bc \neq 0$,
$ u,v \in  \mathfrak b $. We calculate \[
[x,y] =  [e_0, av-cu] + (ad -bc)e_2    \]
which is non zero since the first summand on the right hand side is a multiple of $e_1$. This contradicts
the fact that $\aa \cong \RR ^2$ and part (i) of the proposition follows.

(ii) If $0\to \RR^2 \to \ggo \to  \aff(\RR) \to 0$ splits, then there is
a subalgebra $\hh$ of $\ggo$ isomorphic to $\aff (\RR)$
such that $\ggo=\hh\ltimes \RR^2$.
Set \[ \rho : \hh \to \mathfrak g \mathfrak l (2,\RR), \qquad
\rho(u)=\ad (u)|_{\RR^2} , \quad u\in \hh.\] Then $\rho $ is a Lie algebra
homomorphism.
Let $\mathfrak h =\la x , y \ra, \; [x,y]=y$. If $\rho \equiv 0$ then $\ggo \cong
\aff(\RR) \times \RR^2 =\RR \times \mathfrak r_{3,0}$.  If $\dim $ Im $\rho = 1$,
then $0= [\rho(x), \rho(y)]=\rho([x,y])=\rho(y)$ and $\rho(x)$ is given as follows:
\[ \begin{pmatrix} \mu &0 \\0&\lambda \end{pmatrix},     \quad  \lambda \neq 0,  \qquad
 \begin{pmatrix} \lambda &1 \\0&\lambda\end{pmatrix}   \qquad \text{ or }  \qquad
\begin{pmatrix} \alpha & \beta \\-\beta& \alpha \end{pmatrix}, \quad \beta \neq 0.
 \]
The first possibility gives $\ggo \cong  \RR \times \mathfrak r _{3, \lambda}$ in case
$\mu =0$ and $\ggo \cong \mathfrak r_{4, \mu, \lambda}$ if $\mu \neq 0$. The second possibility
yields $\ggo \cong \rr _{4, \lambda}$ and the last one gives $\ggo \cong \mathfrak r'_{4, 1/\beta, \alpha /\beta}$.

If $\dim $ Im $\rho = 2$, then $\rho(x), \; \rho(y)$ are linearly independent and since
$\ggo '$ is nilpotent and $y \in \ggo '$, we may assume that
\[      \rho(y)=  \begin{pmatrix} 0&1 \\0&0 \end{pmatrix}. \]
It follows from $[\rho(x),\rho(y)]=\rho(y)$ that $\rho(x)$ takes the following
form:
\[  \rho(x)=  \begin{pmatrix} \alpha + 1/2&\beta \\0&\alpha -  1/2 \end{pmatrix}. \]
We can take $\beta=0$ by replacing $x $ with $x-\beta y$.
Let us denote by $\ggo _{\alpha}$ the Lie algebra corresponding to
\[  \rho(x)=  \begin{pmatrix} \alpha + 1/2& 0 \\0&\alpha -  1/2 \end{pmatrix},
\qquad \qquad \rho(y)=  \begin{pmatrix} 0&1 \\0&0 \end{pmatrix}. \]
The following table gives the possibilities for $\ggo _{\alpha}$ according to the
parameter $\alpha$:
\begin{center}
\begin{tabular}{|c|l|} \hline
$\begin{array}{c} \vspace{-.35cm} \\ \alpha  \\  \vspace{-.35cm} \end{array} $ & \quad   $\;\; \ggo_{\alpha}$\\
\hline
 $\begin{array}{c} \vspace{-.35cm} \\ - 1/2 \\  \vspace{-.35cm} \end{array}$ & $\;\;\;\; \dd _4$ \\
$ \begin{array}{c} \vspace{-.35cm} \\  1/2 \\  \vspace{-.35cm} \end{array}$ & $\;\;\;\; \dd _{4, 1}$ \\
$   \alpha \in (-1/2 , 1/2)\cup  (1/2 , 3/2]    $ &
$\;\; \begin{array}{c} \vspace{-.35cm} \\ \dd_{4, \lambda}, \; \lambda= \dfrac 2{2\alpha +1} \\  \vspace{-.35cm} \end{array} $\\
$  \alpha \in (-\infty , - 1/2)\cup (  3/2, \infty )  $ &
$\;\; \begin{array}{c} \vspace{-.35cm} \\ \dd_{4, \lambda}, \; \lambda= \dfrac {\alpha - 1/2}{\alpha +1/2}  \\  \vspace{-.35cm} \end{array} $ \\
\hline
\end{tabular}
\end{center}

This completes the proof of the proposition.
\end{proof}

\subsection{Semidirect extensions of $\aff(\RR)$}

Assume that $\ggo$ contains $\aff(\RR)$ as an ideal and that the short exact sequence
\[ 0\to \aff(\RR) \to \ggo \to \ggo / \aff(\RR)\to 0 \]
splits. The next result states that $\ggo$ is a direct product, that is, $\ggo$ is
isomorphic to $\RR^2 \times \aff(\RR)$ or $ \aff(\RR) \times \aff(\RR)$.
The precise statement is the following:

\begin{prop} Let $\ggo$ be a four dimensional solvable Lie algebra.

(i) If there is a split exact sequence \[ 0\to \aff(\RR) \to \ggo \to \RR^2 \to 0\]
then $\ggo \cong \RR \times \rr _{3,0}= \RR^2 \times \aff(\RR)$.

(ii) If there is a split exact sequence \[ 0\to \aff(\RR) \to \ggo \to \aff(\RR) \to 0\]
then $\ggo \cong \aff (\RR) \times \aff(\RR)$.

\end{prop}
\begin{proof}
Let $\aff(\RR)= \la z,w \ra$ with $[z,w]=w$, then the algebra of derivations
is given as follows:
\[ \text{Der } \aff(\RR)= \left\{ \begin{pmatrix} 0&0\\
a&b\end{pmatrix}, \; a, b \in \RR \right \} \]
with respect to $\{z,w\}$.
If $0\to \aff(\RR) \to \ggo \to \ggo / \aff(\RR) \to 0$ splits, then there is
a subalgebra $\hh$ of $\ggo$ isomorphic to $\ggo / \aff(\RR)$
such that $\ggo=\hh\ltimes \aff(\RR)$.
Set \[ \rho : \hh \to \text{Der } \aff(\RR), \qquad
\rho(u)=\ad (u)|_{\aff(\RR)} , \quad u\in \hh.\] Then $\rho $ is a Lie algebra
homomorphism.

(i)  In this case $\hh= \RR^2$, so the image of $\rho$ is an abelian subalgebra
of $\text{Der } \aff(\RR)$, hence it is one dimensional. Let $\RR^2 = \la x, y \ra$, $\aff(\RR)= \la z,w \ra$.
We may assume that $\rho(y)=0$. Let $\rho(x)=\begin{pmatrix}0&0 \\ a&b \end{pmatrix}$.
If $a=b=0$ the assertion follows. If $b \neq 0$ we may assume that $b=1$
and we can reduce to $a=0$ by changing $z$ to
$z- aw$. Hence, we may assume that the only non zero brackets are $[x,w]=w,
\; [z,w]=w$ and therefore $\ggo =\la x-z , y\ra \times \la z, w \ra$
 where $\la x-z , y\ra \cong \RR^2$. If $b=0, \; a\neq 0$, we may assume that $a=1$,
therefore $\ggo =\la x+w , y\ra \times \la z, w \ra$
 where $\la x+w , y\ra \cong \RR^2$ and the desired assertion follows.

(ii) We have $\hh= \aff(\RR)=\la x, y \ra$, $[x,y]=y$, and the following possibilities for $\rho$:
\[    \rho(x)=\begin{pmatrix} 0&0 \\a&1\end{pmatrix},     \quad
\rho(y)=  \begin{pmatrix} 0&0 \\1&0\end{pmatrix} \qquad
\text{ or } \qquad \rho(x)=\begin{pmatrix} 0&0 \\a&b\end{pmatrix},     \quad
\rho(y)=  0.      \]
We show next that, in both cases,  $\ggo \cong \aff (\RR ) \times \aff(\RR)$.

If the first possibility occurs, take $\la x-z+aw , y+w \ra$ and $ \la  z-aw, w \ra$.
These are complementary ideals isomorphic to $\aff(\RR)$, hence $\ggo \cong
\aff(\RR) \times \aff(\RR)$.

In the second case, take $\la x-bz +aw, y \ra$ and $\la z, w \ra$, which are
ideals isomorphic to $\aff(\RR)$, and  the desired assertion follows.
\end{proof}

\subsection{Product structures of type $\RR \bowtie \hh$}

\
We exhibit in Table~\ref{1,3}  realizations of the Lie algebras obtained in Theorem~\ref{classes}
 as double Lie algebras where $\hh$ is a three dimensional  subalgebra. Note that
the problem of finding such a decomposition is equivalent to the determination
of the three dimensional subalgebras.

\begin{table}
\begin{center}
\begin{tabular}{|l|l|}\hline
${\begin{array}{c} \vspace{-.35cm}\\ \;\;\RR \bowtie \hh \\
\vspace{-.35cm}
\end{array}}$ & $\;\;\;\;\ggo $    \\
\hline \hline ${\begin{array}{c} \vspace{-.35cm}\\ \RR \bowtie \RR^3 \\
\vspace{-.35cm}
\end{array} }$ & $\RR^4, \; \RR \times \hh_3, \; \RR \times \rr_3, \; \RR \times \rr_{3,
\lambda}, \;\RR \times \rr'_{3, \lambda}, \; \nn_4 , \; \rr _4, \; \rr_{4, \lambda}, \;
 \rr_{4, \mu ,\lambda}, \; \rr'_{4, \mu ,\lambda}$\\
\hline
${\begin{array}{c} \vspace{-.35cm}\\ \RR \bowtie \hh_3
\\ \vspace{-.35cm}
\end{array} }$ &
$\RR \times \mathfrak h_3 \, , \; \nn_4 , \; \rr_{4,0}, \; \mathfrak d_4, \;\; \mathfrak d_{4, \lambda}$ ({\scriptsize $
\lambda \geq 1/2 $}), $\;\; \mathfrak d'_{4, \lambda}$ ({\scriptsize $\lambda \geq 0$}),
$\;\; \mathfrak h_4$

 \\
\hline
${\begin{array}{c} \vspace{-.35cm}\\ \RR \bowtie \rr_3
\\ \vspace{-.35cm}
\end{array} }$ &
$ \RR \times \rr_3, \; \rr_4, \; \rr_{4, \lambda}$, ({\scriptsize $\lambda \neq 0$}), $\dd_{4,1}$

\\
\hline
${\begin{array}{c} \vspace{-.35cm}\\  \RR \bowtie \mathfrak r _{3,0} \\
\vspace{-.35cm}
\end{array} }$ &
$\RR \times \mathfrak r _{3,\lambda}, \; \aff(\RR)\times\aff(\RR), \;\;  \mathfrak r _{4,0}, \; \dd_4, \; \dd_{4,1}$\\
\hline
${\begin{array}{c} \vspace{-.35cm}\\ \RR \bowtie \rr'_{3,0} \\
 \vspace{-.35cm}
\end{array} }$ &
$\RR \times \rr'_{3,0}, \; \aff(\CC)$\\

\hline
${\begin{array}{c} \vspace{-.35cm}\\ \RR \bowtie \mathfrak r _{3,\lambda} \\
\vspace{-.35cm}
\end{array} }$ & $\RR \times \rr_{3,\lambda}, \; \aff(\RR)\times \aff(\RR)$
 ({\scriptsize $\lambda =-1$}), $\aff(\CC)$ ({\scriptsize $\lambda =1$}),
   $\; \rr_{4,\lambda},\; \rr_{4 , \mu , \lambda}, \; \hh_4$  ({\scriptsize $\lambda =2$}),
$ \dd _{4, \lambda}, \; \dd_{4,1-\lambda}$

  \\
\hline
${\begin{array}{c} \vspace{-.35cm}\\ \RR \bowtie \mathfrak r' _{3,\lambda} \\
\vspace{-.35cm}
\end{array} }$ & $\RR \times \rr'_{3,\lambda}, \; \aff(\CC), \; \rr' _{4, \mu ,\lambda}$
  \\
\hline
\end{tabular}
\bigskip
\end{center} \caption{} \label{1,3}
\end{table}

\begin{itemize}

\item $\RR \bowtie \RR^3$:
The Lie algebras $\RR \times \hh_3  , \;
\RR \times \rr_3, \; \RR \times \rr_{3,
\lambda}, \;\RR \times \rr'_{3, \lambda}, \; \nn_4 , \; \rr _4, \;
\rr_{4, \lambda}, \;
 \rr_{4, \mu ,\lambda}$ and $ \rr'_{4, \mu ,\lambda}$ were obtained
in Theorem~\ref{classes} as semidirect extensions of $\RR^3$.

\item $\RR \bowtie \hh_3$:  The Lie algebras $\RR \times \mathfrak h_3 \, , \; \nn_4 , \;  \; \mathfrak d_4, \;\; \mathfrak d_{4, \lambda}$ ({\scriptsize $
\lambda \geq 1/2 $}), $\;\; \mathfrak d'_{4, \lambda}$ ({\scriptsize $\lambda \geq 0$}),
 and $\;\; \mathfrak h_4$ were obtained in Theorem~\ref{classes} as semidirect extensions
 of $\hh_3$. On the other hand, \[\rr_{4,0} \cong \la e_1 \ra \bowtie \la e_0, e_2, e_3 \ra.\]

\item $\RR \bowtie \rr_3$:
\begin{eqnarray*} \rr_4 &\cong & \la e_3 \ra \bowtie \la e_0, e_1 , e_2 \ra ,\\
 \rr_{4,\lambda} &\cong & \la e_1 \ra \bowtie \la e_0, e_2 , e_3 \ra, \qquad \lambda \neq 0,\\
 \dd_{4,1}& \cong & \la e_0 \ra \ltimes \la  e_0+ e_2, e_1, e_3 \ra . \end{eqnarray*}

\item $\RR \bowtie \rr_{3,0}$: \begin{eqnarray*} \aff(\RR)\times \aff(\RR) & \cong &
\la e_0 \ra \ltimes \la e_1, e_2 , e_3\ra , \\
\rr_{4,0} &\cong & \la e_3 \ra \ltimes \la e_0, e_1, e_2 \ra, \\
\dd_{4} & \cong & \la e_2 \ra \bowtie \la e_0 , e_1,  e_3\ra, \\ \dd_{4,1} &\cong &
\la e_1 \ra \bowtie \la e_0, e_2 , e_3\ra.  \end{eqnarray*}

\

\item $\RR \bowtie \rr'_{3,0}$: $\aff (\CC) \cong \la e_0 \ra \ltimes \la e_1, e_2 , e_3\ra$.

\

\item $\RR \bowtie \mathfrak r _{3,\lambda}$: $\aff(\RR) \times \aff(\RR)$ was obtained
in Theorem~\ref{classes} as a semidirect extension of $\rr_{3, -1}$.
  \begin{eqnarray*}
\hh_4&=& \la e_2 \ra \bowtie \la e_0, e_1, e_3  \ra \cong \RR \bowtie \rr_{3,2}, \qquad
\lambda=2,\\
 \aff(\CC)&=& \la e_1\ra \ltimes   \la e_0, e_2, e_3  \ra  \cong
\RR \ltimes \rr_{3,1}, \qquad \lambda =1, \\
\rr_{4, \lambda} &\cong &\la e_3 \ra \bowtie \la e_0, e_1 ,e_2 \ra,\\
 \rr_{4,\mu , \lambda}
&\cong & \la e_2 \ra \bowtie \la e_0, e_1, e_3 \ra , \\
\rr_{4,\mu , \lambda}
& \cong & \RR \bowtie \rr_{4, \mu}= \la e_3 \ra \bowtie \la e_0, e_1, e_2 \ra, \\
 \dd _{4, \lambda} & \cong & \la e_2 \ra \bowtie \la e_0, e_1 , e_3 \ra , \; \dd_{4,1-\lambda} \cong \la e_1 \ra \bowtie \la e_0, e_2 , e_3 \ra. \end{eqnarray*}

\item $\RR \bowtie \mathfrak r '_{3,\lambda}$: \begin{eqnarray*} \aff(\CC) &\cong &\la e_0 \ra \ltimes
\la \lambda e_0 -e_1 , e_2, e_3 \ra, \\ \rr' _{4, \mu ,\lambda} &\cong &\la e_1 \ra \bowtie \la e_0, e_2 , e_3 \ra . \end{eqnarray*}
\end{itemize}

\

\section{Applications: Manin triples and complex product  structures}

\subsection{Manin triples on 4-dimensional solvable Lie algebras}

An important example of double Lie algebras are Manin triples \cite{LW}. We recall that a {\it Manin triple} is a double Lie algebra $(\ggo, \ggo_+,\ggo_-)$ with an invariant metric, that is, a non degenerate symmetric bilinear form $(\; ,\, )$ which satisfies:
\[ ([x,y],z) + (y,[x,z]) = 0 \qquad \text{ for all } x, y, z \in \ggo\]
such that $\ggo_+$ and $\ggo_-$ are isotropic subalgebras. In particular $\ggo=\ggo_+ \bowtie \ggo_-$, where $\ggo_+$ and $\ggo_-$ have the same dimension. Thus, Manin triples are special cases of paracomplex structures.

The next proposition makes use of \cite{B-K} and the results of the previous section to
obtain that there is only one four dimensional solvable Lie algebra giving rise to Manin triples.

\begin{prop} Let $(\ggo,\ggo_+, \ggo_-)$ be a Manin triple such that $\ggo$ is a
non abelian four dimensional solvable Lie algebra. Then $\ggo$ is
isomorphic to $\dd_4$ with the invariant metric given by:
\[ \alpha = (e_0, e_3) = (e_1,e_2), \quad \alpha \ne 0,\]
where the isotropic subalgebras $\ggo_+$ and $\ggo_-$ are given as
follows:

$\, $(i) $\, \ggo_+ = \la e_0+\mu e_2, e_1-\mu e_3\ra , \quad
\ggo_- = \la e_0 +\nu e_2, e_1-\nu e_3\ra$ with $\mu \ne \nu$; or

(ii)  $\ggo_+= \la e_0+\mu e_2, e_1-\mu e_3\ra , \quad \ggo_- =
\la e_2, e_3\ra$.
\end{prop}
\begin{proof} According to \cite{B-K} a non abelian solvable Lie algebra which admits
an invariant metric is isomorphic either to $\dd'_{4,0}$ or
$\dd_4$. It was proved in Theorem~\ref{paracomplex} that
$\dd'_{4,0}$ does not admit paracomplex structures. Thus, we need
to investigate the possible paracomplex structures on $\dd_4$.
It is easy to see that the metric on $\dd_4$  given by:
$$(e_0,e_3)=(e_1,e_2)=\alpha, \quad \text{ with } \alpha \ne 0,$$
is invariant.
Any two-dimensional  isotropic non abelian subalgebra of $\dd_4$
  is isometrically isomorphic to:
\[\la e_0+\mu e_2, e_1-\mu e_3\ra,\]
where the isometric isomorphism is given by $\phi(e_1)=e_2, \; \phi (e_i)=-e_i, \; i=0,3$.
On the other hand, any two-dimensional isotropic abelian subalgebra is
isometrically isomorphic to:
\[\la e_2, e_3\ra.\]
It follows from \ref{no22} that $\dd_4$ does not admit a
decomposition of type $\RR^2 \bowtie \RR^2$.  If both $\ggo_+$ and
$\ggo_-$ are isomorphic to $\aff(\RR)$ then we are led to case
(i). In case $\ggo _+ \cong \aff(\RR)$ and $\ggo_- \cong \RR^2$ we
obtain case (ii), and the proposition follows.
\end{proof}

\smallskip

\subsection{Complex product structures on four dimensional solvable Lie algebras}

In this subsection we determine all four dimensional solvable Lie
algebras which admit a complex product structure (see Table \ref{cps}), using the classification of complex structures on this class of Lie algebras given in \cite{SJ,O1} together with the results in \S 2.2. We give in this way an alternative proof of a result by Blazi\'c and Vukmirovi\'c (\cite{BV}), where complex product structures were referred to as {\em para-hypercomplex structures}.

We recall that a complex structure on a Lie algebra $\ggo$ is an endomorphism $J:\ggo\rightarrow\ggo$ such that $J^2=-\,$Id and
\[J[x,y]=[Jx,y]+[x,Jy]+J[Jx,Jy]\]
for all $x,y\in\ggo$. A {\em complex product structure} on a Lie algebra
$\ggo$ is a pair $\{J,E\}$ where $J$ is a complex structure and
$E$ is a product structure on $\ggo$ such that $JE=-EJ$. This is
equivalent to having a splitting of $\ggo$ as
$\ggo=\ggo_+\oplus\ggo_-$, where $\ggo_+$ and $\ggo_-$ are Lie
subalgebras of $\ggo$ such that $\ggo_-=J\ggo_+$. From this it
follows that $E$ is, in fact, a paracomplex structure on $\ggo$.

At this point we refer the reader to Table \ref{cps}.

\begin{table}
\begin{center}
\begin{tabular}{|c|c|c|}\hline

$\begin{array}{c} \vspace{-.35cm}\\ \text{Lie algebra} \\ \vspace{-.35cm} \end{array}$ & Complex structure & Paracomplex structure \\ \hline

$\begin{array}{c} \vspace{-.35cm}\\ \aff(\RR)\times\aff(\RR) \\ \vspace{-.35cm} \end{array}$ & $Je_0=e_3,\,Je_1=e_2$ & $\la e_0,e_1\ra\ltimes\la e_2,e_3\ra$ \\
\hline

$\begin{array}{c} \vspace{-.35cm}\\ \RR\times\hh_3 \\ \vspace{-.35cm} \end{array}$ & $Je_0=-e_3,\,Je_1=e_2$ & $\la e_0,e_1\ra\ltimes\la e_2,e_3\ra$ \\
\hline

$\begin{array}{c} \vspace{-.35cm}\\ \RR\times\rr_{3,0} \\ \vspace{-.35cm} \end{array}$ & $Je_3=e_0,\,Je_1=e_2$ & $\la e_1,e_3\ra\ltimes\la e_0,e_2\ra$ \\
\hline

$\begin{array}{c} \vspace{-.35cm}\\ \RR\times\rr_{3,1} \\ \vspace{-.35cm} \end{array}$ & $Je_0=e_1,\,Je_2=e_3$ & $\la e_1,e_3\ra\ltimes\la e_0,e_2\ra$ \\
\hline

$\begin{array}{c} \vspace{-.35cm}\\ \aff(\CC) \\ \vspace{-.35cm} \end{array}$ & $Je_0=e_2,\,Je_2=e_3$ & $\la e_0,e_1\ra\ltimes\la e_2,e_3\ra$ \\
\hline

$\begin{array}{c} \vspace{-.35cm}\\ \rr_{4,1} \\ \vspace{-.35cm} \end{array}$ & $Je_0=e_3,\,Je_1=e_2$ & $\la e_0,e_1\ra\ltimes\la e_2,e_3\ra$\\
\hline

$\begin{array}{c} \vspace{-.35cm}\\ \rr_{4,\lambda,\lambda} \\ \vspace{-.35cm} \end{array}$, \scriptsize{$\lambda\neq 0$} & $Je_0=e_1,\,Je_2=e_3$ &
$\la e_0,e_2\ra\ltimes\la e_1,e_3\ra$ \\
\hline

$\begin{array}{c} \vspace{-.35cm}\\ \rr_{4,\mu,1} \\ \vspace{-.35cm} \end{array}$, \scriptsize{$\mu\neq 0,\,\pm1$} &
$Je_0=e_2,\,Je_1=e_3$ & $\la e_0,e_1\ra\ltimes\la e_2,e_3\ra$ \\
\hline

$\begin{array}{c} \vspace{-.35cm}\\ \rr'_{4,\mu,\lambda} \\ \vspace{-.35cm} \end{array}$ & $Je_0=e_1,\,Je_2=e_3$ & $\la e_0,e_1\ra\ltimes\la e_2,e_3\ra$ \\
\cline{2-3}            & $\begin{array}{c} \vspace{-.35cm}\\ Je_0=e_1,\,Je_2=-e_3 \\ \vspace{-.35cm} \end{array}$ & $\la e_0,e_1\ra\ltimes\la e_2,e_3\ra$\\
\hline

$\begin{array}{c} \vspace{-.35cm}\\ \dd_4 \\ \vspace{-.35cm} \end{array}$ & $Je_0=-e_1,\,Je_2=e_3$ & $\la e_0,e_1\ra\ltimes\la
e_2,e_3\ra$ \\ \cline{2-3}
& $\begin{array}{c} \vspace{-.35cm}\\ Je_0=e_3-e_1,\,Je_1=e_0-e_2,\,Je_2=e_3 \\ \vspace{-.35cm} \end{array}$ & $\la e_0,e_2\ra\ltimes\la e_1,e_3\ra$ \\
\hline

$\begin{array}{c} \vspace{-.35cm}\\ \dd_{4,1} \\ \vspace{-.35cm} \end{array}$ & $Je_0=e_1,\,Je_2=-e_3$ & $\la e_0,e_2\ra\ltimes\la e_1,e_3\ra$ \\
\hline

$\begin{array}{c} \vspace{-.35cm}\\ \dd_{4,1/2} \\ \vspace{-.35cm} \end{array}$ & $Je_0=e_3,\,Je_1=e_2$ & $\la e_0,e_1\ra\ltimes\la e_2,e_3\ra$ \\ \cline{2-3}
              & $\begin{array}{c} \vspace{-.35cm}\\ Je_0=e_3,\,Je_1=-e_2 \\ \vspace{-.35cm}\end{array}$ & $\la e_0,e_1\ra\ltimes\la e_2,e_3\ra$ \\ \cline{2-3}
              & $\begin{array}{c} \vspace{-.35cm}\\ Je_0=e_1,\,Je_2=-2e_3 \\ \vspace{-.35cm} \end{array}$ & $\la e_0,e_2\ra\ltimes\la e_1,e_3\ra $\\ \hline

$\begin{array}{c} \vspace{-.35cm}\\ \dd_{4,\lambda} \\ \vspace{-.35cm} \end{array}$, \scriptsize{$\lambda\neq 1,1/2$} &
$Je_0=(1-\lambda)e_2,\,Je_1=e_3 $ & $\la e_0,e_1\ra\ltimes\la
e_2,e_3\ra$ \\ \cline{2-3}
& $\begin{array}{c} \vspace{-.35cm}\\ Je_0=-\lambda e_1,\,Je_2=e_3 \\ \vspace{-.35cm}\end{array}$ & $\la e_0,e_2\ra\ltimes\la e_1,e_3\ra$ \\
\hline

$\begin{array}{c} \vspace{-.35cm}\\ \hh_4 \\ \vspace{-.35cm} \end{array}$ & $Je_0=4e_2,\,Je_1=4e_3$ & $\la e_0,e_1\ra\bowtie\la e_2,e_3\ra$ \\ \hline
\end{tabular}
\bigskip
\caption{Complex product structures}\label{cps}
\end{center}
\end{table}

\

\begin{remarks}
\item[a)] The Lie algebras $\RR\times\rr'_{3,\lambda}$ admit
complex structures (see \cite{SJ}) and paracomplex structures (see
Table~\ref{pc}). Nevertheless, they do not admit any complex
product structure. To show this, we state the following result,
which is proved in \cite{AD}:

\begin{prop}
Let $\{J,E\}$ be a complex product structure on the Lie algebra
$\ggo$ and let $(\ggo,\ggo_+,\ggo_-)$ be the associated double Lie
algebra. Then the following assertions are equivalent:
\begin{enumerate}
\item[$(i)$] $J$ is an abelian complex structure, i.e., $[Jx,Jy]=[x,y]$ for all $x,y\in\ggo$.
\item[$(ii)$] The Lie subalgebras $\ggo_+$ and $\ggo_-$ are abelian;
\item[$(iii)$] $E$ is an abelian product structure, i.e., $[Ex,Ey]=-[x,y]$ for all $x,y\in\ggo$.
\end{enumerate}
\end{prop}

It is known that the Lie algebra $\RR\times\rr'_{3,\lambda}$ does not
admit any abelian complex structure (see \cite{SJ}). However,
from Lemma \ref{sub-aff} this Lie algebra admits only abelian paracomplex
structures and thus, from the previous proposition, there is no complex product
structure on $\RR\times\rr'_{3,\lambda}$.

\item[b)] The Lie algebra $\aff(\CC)$ admits other complex structures, given by:
\[ J_{\alpha,\beta}e_0=\frac{\alpha}{\beta}e_0+\frac{\alpha^2+\beta^2}{\beta}e_1,\quad J_{\alpha,\beta}e_2=e_3 \]
with $\alpha\in\RR,\,\beta\in\RR\setminus\{0\}$. However, there is no
paracomplex structure on $\aff(\CC)$ which anticommutes with $J_{\alpha,\beta}$. Let
us show this last assertion. It is known from Proposition \ref{nose} that $\aff(\CC)$
does not admit decompositions of type $\aff(\RR)\bowtie\aff(\RR)$. Also, since
the complex structure $J_{\alpha,\beta}$ is not abelian, any complex product
structure on $\aff(\CC)$ induces a decomposition of type $\aff(\RR)\bowtie\RR^2$. Let $\hh$
be a subalgebra of $\aff(\CC)$ isomorphic to $\aff(\RR)$; then $\hh$ has a basis
$u=e_0+a_2e_2,\,v=b_2e_2+b_3e_3$. Then $J_{\alpha,\beta}u=\frac{\alpha}{\beta}e_0+\frac{\alpha^2+\beta^2}{\beta}e_1-a_3e_2+a_2e_3,\, J_{\alpha,\beta}v=-b_3e_2+b_2e_3$. Then we must have $[Ju,Jv]=0$ and from this we obtain the system
\[ \begin{cases}
     \frac{\alpha^2+\beta^2}{\beta}b_2+\frac{\alpha}{\beta}b_3=0 \\
     \frac{\alpha}{\beta}b_2-\frac{\alpha^2+\beta^2}{\beta}b_3=0
   \end{cases} \]
It is easy to see that the only solution of this system is $b_2=b_3=0$, i.e., $v=0$, a contradiction. Therefore, there are no product structures on $\aff(\CC)$ which anticommute with $J_{\alpha,\beta}$.

\item[c)] The Lie algebras $\dd'_{4,\lambda}$ admit complex structures (see \cite{O1}) but, according to Theorem \ref{paracomplex}, they do not admit any paracomplex structure. Hence, they do not carry complex product structures.

\item[d)] Table \ref{cps} shows examples of complex product structures on $\RR\times\hh$ and $\rr_{4,1,1}$. On the other hand, all equivalence classes of complex product structures on these Lie algebras were determined in \cite{AS}, section 6.2.
\end{remarks}

\

\section*{Appendix I - Matrix realizations}

We exhibit below matrix realizations of the indecomposable Lie algebras listed in Theorem~\ref{classes}, where
indecomposable means that they do not split as a direct product of
lower dimensional Lie algebras.
All matrices have real coefficients.

\begin{center}

\begin{tabular}{rcrc}

{$\mathfrak n_4$:} & {$\left(

\begin{matrix}

0 & x & 0 & w \\

0 & 0 & x & y \\

0 & 0 & 0 & z \\

0 & 0 & 0 & 0

\end{matrix}

\right)$} &

{ $\aff(\CC)$:} & {$\left(

\begin{matrix}

x & z & y  \\

-z & x & w  \\

0 & 0 & 0

\end{matrix}

\right)\,$}

\\

 &&&\\

{ $\mathfrak r _{4}$:} & {$\left(

\begin{matrix}

x& x & 0& y  \\

0 &   x  & x & z  \\

0 & 0 &  x  &  w\\

0 & 0& 0 & 0

\end{matrix}

\right)\,$}&

{ $\mathfrak r _{4,\lambda}$:} & {$\left(

\begin{matrix}

x & 0 & 0 & y \\

0 & \lambda x  & x  & z \\

0 & 0 & \lambda x  & w \\

0 & 0 & 0 & 0

\end{matrix}

\right)$} \\

&&&\\

$\begin{array}{c}

\mathfrak r _{4, \mu ,\lambda}:  \\

\text{ \scriptsize{$\mu \lambda \ne 0$}}, \text{ \scriptsize{$-1 < \mu \leq \lambda \leq 1$ }} \\

\text{or \scriptsize{$-1 = \mu \leq \lambda < 0$ }} \\

\end{array}$

 & {$\left(

\begin{matrix}

x& 0 & 0& y  \\

0 &  \mu x  & 0 & z  \\

0 & 0 & \lambda x  &  w \\

0 & 0& 0 & 0

\end{matrix}

\right)\,$}&

$\begin{array}{c}

\mathfrak r _{4, \mu ,\lambda}': \\

\text{\scriptsize{$\mu > 0 $}}

\end{array}$

& {$\left(

\begin{matrix}

\mu x & 0 & 0 & y \\

0 &  \lambda x &  x& z \\

0 & - x &  \lambda x& w \\

0 & 0 & 0 & 0

\end{matrix}

\right)$} \\

&&&\\

{ $ \mathfrak d _{4} $:} & $\left\{ \quad  \begin{array}{c}
\begin{pmatrix} 0&x &z \\ 0&w&y\\ 0&0&0 \end{pmatrix}
 \\  \\
{\left(

\begin{matrix}

w & 0 & 0 & x \\

0 & -w & 0 & y \\

-\frac12 y & \frac12 x & 0 & z \\

0 & 0 & 0 & 0

\end{matrix}

\right)} \end{array} \right.$ &  $\begin{array}{c} \mathfrak d _{4, \lambda} :\\

  \text{\scriptsize{$\lambda \geq \frac12 $}}

\end{array}$

 &  $\left\{ \quad \begin{array}{c} \left(

\begin{matrix}

w & x & z  \\

0 & (1 - \lambda)w & y \\

0 & 0 & 0

\end{matrix}

\right)  \\    \\

\left(

\begin{matrix}

\lambda w & 0 & 0 & x \\

0 & (1 - \lambda)w & 0 & y \\

-\frac12 y & \frac12 x & w & z \\

0 & 0 & 0 & 0

\end{matrix}

\right)
\end{array} \right.$ \\

&&& \\

$\begin{array}{c} \mathfrak d _{4, \lambda}' :\\

  \text{\scriptsize{$\lambda \geq 0$}}

\end{array}$

 & {$\left(

\begin{matrix}

  \lambda w&   w& 0 & x \\

 -  w &  \lambda w& 0 & y \\

-\frac12 y & \frac12 x & 2\lambda w & z \\

0 & 0 & 0 & 0

\end{matrix}

\right)$}

 & { $\mathfrak h _4$:} & {$\left(

\begin{matrix}

\frac12 w & w & 0 & x \\

0 & \frac12 w & 0 & y \\

-\frac12 y & \frac12 x & w & z \\

0 & 0 & 0 & 0

\end{matrix}\right)$} \\

&&&\\

\end{tabular}

\end{center}

\

\section*{Appendix II - Comparison with previous classifications}

In this section we carry out a comparison with various results  which can be found in the
literature. Our main goal is to establish a correspondence between
the description obtained by other authors and the Lie algebras
appearing in Theorem~\ref{classes}.

\

\subsection{} We start by comparing our results with the ones
obtained by Dozias as appearing in  \cite{Ve}, Table 1.1, p. 180.

\

\begin{center}
\small{
\begin{tabular}{|c|c|c|c|c|c|c| c|c|c|c|c|}\hline
$\begin{array}{c} \vspace{-.2cm} \\ \ggo_{4,1} \\ \vspace{-.2cm} \end{array}$
& $\ggo_{4,2}$ & $\ggo_{4,3}$ & $\ggo_{4,4}$ &
$\ggo_{4,5}(\alpha,\beta)$ & $\ggo_{4,6}(\alpha)$ &
$\ggo_{4,7}$ & $\ggo_{4,8}(\alpha,\beta)$ &$\ggo_{4,9}(0)$ & $\ggo_{4,9}(\alpha)\, , \,\;$
\scriptsize{$\alpha \neq 0$} & $\ggo_{4,10}$ & $\ggo_{4,11}(\alpha)$  \\
\hline
$\begin{array}{c} \vspace{-.2cm} \\ \mathfrak d _{4,0} \\ \vspace{-.2cm} \end{array} $
& $\aff(\CC)$ & $\mathfrak n_{4}$ & $\mathfrak r_{4,0}$ &
$\mathfrak r_{4,\alpha,\beta}$ & $\mathfrak r_{4,\alpha}$ & $\mathfrak r_{4}$ &
$\mathfrak r_{4,\alpha,\beta}'$ &$\mathfrak d_4$ & $\mathfrak d_{4,1-1/{\alpha}}$ &
$\mathfrak h_{4}$ & $\mathfrak d _{4,{\alpha}}'$ \\
\hline
\end{tabular} }
\end{center}

\

\

\subsection{}
We recall  below the classification given by Mubarakzyanov \cite{Mu} and
then we establish the correspondence with the algebras appearing
in Theorem~\ref{classes}.

\

\begin{center}

\small{
\begin{tabular}{|c|llll|}\hline
$\begin{array}{c} \vspace{-.3cm} \\ \text{Notation  in \cite{Mu}} \\ \vspace{-.3cm} \end{array} $
&   & \hspace{2.1cm} Lie bracket  & relations &\\\hline
$\begin{array}{c} \vspace{-.3cm} \\ \mathfrak{g}_{4,1}  \\ \vspace{-.3cm} \end{array} $
& $[e_2 , e_4]= e_1 $& $[e_3,e_4]=e_2$& & \\
$\begin{array}{c} \vspace{-.3cm} \\ \mathfrak{g}_{4,2} \\ \vspace{-.3cm} \end{array}$
& $[e_1 , e_4]= \alpha e_1 $& $[e_2,e_4]=e_2$& $[e_3,e_4]=e_2 +e_3$& \\
$\begin{array}{c} \vspace{-.3cm} \\ \mathfrak{g}_{4,3} \\ \vspace{-.3cm} \end{array}$&
$[e_1 , e_4]=  e_1 $& $[e_3,e_4]=e_2$& & \\
$\begin{array}{c} \vspace{-.3cm} \\ \mathfrak{g}_{4,4} \\ \vspace{-.3cm} \end{array}$
& $[e_1 , e_4]=  e_1 $& $[e_2,e_4]=e_1 + e_2$& $[e_3,e_4]=e_2 +e_3$& \\
$\begin{array}{c} \vspace{-.3cm} \\ \mathfrak{g}_{4,5} \\ \vspace{-.3cm} \end{array}$
& $[e_1 , e_4]= e_1 $& $[e_2,e_4]=\beta e_2$& $[e_3,e_4]=\gamma e_3$& -1$\le \gamma \le \beta \le$ 1, \, $\gamma \beta \ne 0$\\
$\begin{array}{c} \vspace{-.3cm} \\ \mathfrak{g}_{4,6} \\ \vspace{-.3cm} \end{array}$
& $[e_1 , e_4]= \alpha e_1 $& $[e_2,e_4]=pe_2-e_3$& $[e_3,e_4]=e_2 +pe_3$& $\alpha \ne$ 0, p$\ge 0$\\
$\begin{array}{c} \vspace{-.3cm} \\ \mathfrak{g}_{4,7} \\ \vspace{-.3cm} \end{array}$
& $[e_2 , e_3]=  e_1 $& $[e_1,e_4]=2e_1$& $[e_2,e_4]=e_2 $& $[e_3,e_4]=e_2 +e_3$\\
$\begin{array}{c} \vspace{-.3cm} \\ \mathfrak{g}_{4,8} \\ \vspace{-.3cm} \end{array}$
& $[e_2 , e_3]=  e_1 $& $[e_1,e_4]=(1+h)e_1$& $[e_2,e_4]=e_2$& $[e_3,e_4]=h e_3$,$|\text{h}|\le$1\\
$\begin{array}{c} \vspace{-.3cm} \\ \mathfrak{g}_{4,9} \\ \vspace{-.3cm} \end{array}$
& $[e_2 , e_3]=  e_1 $& $[e_1,e_4]=2pe_1$& $[e_2,e_4]=pe_2 -e_3$& $[e_3,e_4]=e_2 +pe_3$, p$\ge$0\\
$\begin{array}{c} \vspace{-.3cm} \\ \mathfrak{g}_{4,10} \\ \vspace{-.2cm} \end{array}$
& $[e_1 , e_3]=  e_1 $& $[e_2,e_3]=e_2$& $[e_1,e_4]=-e_2$& $[e_2,e_4]=e_1 +e_3$\\\hline
\end{tabular}}

\end{center}

\

The correspondence is as follows:

\

\begin{center}
\small{
\begin{tabular}{|c| c|c|c|c|c|c| c| c| c|}\hline
$\begin{array}{c} \vspace{-.3cm} \\ \mathfrak g_{4,1} \\ \vspace{-.3cm} \end{array}$
& ${\mathfrak g_{4,2}}$ & $\mathfrak g_{4,3}$ & $\mathfrak g_{4,4}$ &
${\mathfrak g_{4,5}}$ & ${\mathfrak g_{4,6}}$ & $\mathfrak g_{4,7}$ & $\mathfrak g _{4,8}$ &${\mathfrak g_{4,9}}$ & $\mathfrak g_{4,10}$  \\
\hline
$\begin{array}{c} \vspace{-.3cm} \\ \mathfrak n_4 \\ \vspace{-.3cm} \end{array}$
& $\mathfrak r_{4,\alpha}$ & $\mathfrak r_{4,0}$ &
$\mathfrak r_{4}$ & $\mathfrak r_{4,\beta,\gamma}$ & $\mathfrak r_{4,\alpha,p}'$ &
$\mathfrak h_{4}$ & $\mathfrak d _4,\mathfrak d_{4,1/1+b}$ &$ \mathfrak d_{4,a}'$ & $\aff (\CC)$\\
\hline
\end{tabular} }
\end{center}

\

\

\subsection{}
In \cite{Z} invariants of real Lie algebras of dimension at most
five are given. In particular, a list of four dimensional solvable
Lie algebras, based on that of \cite{Mu}, is shown in Table I, p.
988.  The relation with Theorem~\ref{classes} is :

\

\begin{center}
\small{
\begin{tabular}{|c|c|c|c|c|c|c| c|c|c|c|c|}\hline
$\begin{array}{c} \vspace{-.3cm} \\ A_{4,1} \\ \vspace{-.3cm} \end{array}$
& ${A_{4,2}}^a$ & $A_{4,3}$ & $A _{4,4}$ &
${A_{4,5}}^{a, b}$ & ${A_{4,6}}^{a, b}$ & $A _{4,7}$ & $A _{4,8}$  &${A_{4,9}}^b$ & $A_{4,10}$ & ${A_{4,11}}^a$ & $A_{4,12}$ \\
\hline
$\begin{array}{c} \vspace{-.3cm} \\ \mathfrak n_4 \\ \vspace{-.3cm} \end{array}$
& $\mathfrak r_{4,a}$ & $\mathfrak r_{4,0}$ &
$\mathfrak r_{4}$ & $\mathfrak r_{4,a,b}$ & $\mathfrak r_{4,a,b}'$ &
$\mathfrak h_{4}$ & $\mathfrak d _4$& $\mathfrak d_{4,1/1+b}$ & $\mathfrak d_{4,0}'$ &  $\mathfrak d_{4,a}'$ & $\aff (\CC)$\\
\hline
\end{tabular} }
\end{center}

\

\subsection{}
The classification of complex structures on four dimensional Lie
algebras was carried out by Snow in \cite{SJ} and by Ovando in
\cite{O1}. To achieve this classification a description is given
in \cite{SJ}, p. 400, of four dimensional solvable Lie algebras
when the commutator ideal has dimension $1$ or $2$. We compare
below the list given by Snow with the one obtained in
Theorem~\ref{classes}.

\

\begin{center}

\small{

\begin{tabular}{|c| c|c|c|c|c|c| c|}\hline
$\begin{array}{c} \vspace{-.3cm} \\ S1 \\ \vspace{-.3cm} \end{array}$
& $S2$& $\; \;S3 \;\;$ & $\;\;S4\;\;$& $S5 _d \, , \,\;${\scriptsize $d\neq 0$}  & $S6$ &
$S7_{0,c}\, , \,\;$ {\scriptsize $c> 0$}  &  $S7_{1,c}\, , \,\;$ {\scriptsize $4c>1$} \\
\hline
$\begin{array}{c} \vspace{-.3cm} \\ \RR \times \mathfrak h _3 \\ \vspace{-.3cm} \end{array}$
& $\RR^ 2\times \aff (\RR)$  &
$\mathfrak r _{4,0}$ & $\mathfrak n _4$ & $\RR \times \mathfrak r_{3, d}$&
$\RR \times \mathfrak r _3$ & $\RR \times \mathfrak r '_{3,0}$ &
$\RR \times \mathfrak r '_{3,\sqrt{4c-1}}$ \\
\hline
\end{tabular} }
\end{center}

\

\

\begin{center}
\small{
\begin{tabular}{|c|c|c|c|c|}\hline
 $\begin{array}{c} \vspace{-.3cm} \\ \;\;S8\;\; \\ \vspace{-.3cm} \end{array}$
& $\quad S9 \quad $ & $S10 _{d,d}\, , \,\; d\neq 0 $ & $S10_{d,c}$, {\scriptsize $\; c\neq d, \; d\neq 0$} & $S11_{d,c}\, , \;$
{\scriptsize $d^2-4c <0, \; d=0, 1$} \\
 \hline
 $\begin{array}{c} \vspace{-.3cm} \\ \aff(\RR)\times \aff(\RR) \\ \vspace{-.3cm} \end{array}$
& $\mathfrak d_{4,1 }$&
 $\RR \times \mathfrak r_{3, d}$ &  $\aff (\RR)\times \aff (\RR)$ &
 $\aff (\CC)  $\\
 \hline
\end{tabular} }
\end{center}

\

The above correspondence shows that some of the families appearing
in \cite{SJ} become a single Lie algebra. Also, there exist
isomorphisms between different families. We give below the
proof of these statements.

\begin{itemize}

\item {\bf S7}

We recall from \cite{SJ} the definition of the Lie algebra
$S7_{d,c}\, , \,\;$ {\scriptsize $d^2-4c <0, \; d=0$ or $1$} with
basis $x,y,z,w$:
\[ [x,y]=w, \qquad\qquad [x,w]=-cy+dw.\]
Observe that if $d=0$ then $c>0$ and $\ad (x)_{|\ggo '}$ has
eigenvalues $\pm ic$. We can take a real basis of $\ggo '$ such
that $\ad (x)_{|\ggo '}$ takes the form $\begin{pmatrix} \; 0 & c
\\ -c & 0
\end{pmatrix}$.
Changing $x$ by $x/c$ we see that $S7_{0,c} \cong \RR \times
\mathfrak{e}(2)$ for all $c>0$. If $c=1$, then $\ad (x)_{|\ggo '}$
has eigenvalues $1/2 \pm i \lambda /2, \;$ where $\lambda=
\sqrt{4c-1}$. Taking $x'=x/2$, there exists a real basis of $\ggo
'$ such that $\ad (x')_{|\ggo '}$ takes the form $\begin{pmatrix}
\; 1 & \lambda \\ -\lambda & 1 \end{pmatrix}$, hence $S7_{1,c}
\cong \RR \times \mathfrak r '_{3,\lambda}$ for all $c$ such that
$1-4c<0$.

\item {\bf S10}

Consider next the Lie algebra $S{10}_{d,c}\,,\,\;$ {\scriptsize $c,d \in \RR, d\neq 0$}:
\begin{equation} [x,y]=y, \qquad [x,w]=dw, \qquad [z,y]=y,
\qquad [z,w]=cw. \label{s10} \end{equation}
If $c=d$, then changing $z$ by $x-z$, we see that  $S{10}_{d,d}
\cong \RR \times \mathfrak r_{3, d}$ for all $d\neq 0$.

If $c \neq d$, let $x', y', z', w'$ be the basis of $S{10}_{1,0}$
satisfying~\eqref{s10} and $x, y, z, w$ the corresponding basis
of $S{10}_{d,c}$. Define a linear map $\psi: S{10}_{d,c} \to   S{10}_{1,0}$  by
\[  \psi(x)=x' +(d-1)z', \qquad \psi (y)=w', \qquad
\psi(z)=x' +(c-1)z', \qquad \psi(w)=y'.
\]
It turns out that $\psi$ is a Lie algebra isomorphism for all $c\neq d$
and therefore $S{10}_{d,c} \cong S{10}_{1,0} \cong \aff (\RR)
\times \aff(\RR)$, where the last isomorphism follows by changing
$x'$ to $x'-z'$.

\item {\bf S11}

Consider the Lie algebra $S11_{d,c}\, , \;${\scriptsize $d^2-4c <0, \; d=0, 1$}:
\[ [x,y]=y, \qquad\qquad [x,w]=w, \qquad\qquad [z,y]=w,
\qquad\qquad [z,w]=-cy+dw.\] If $d=0$, then $\ad (z/c)_{|\ggo '}$
has eigenvalues $\pm i $ and there exists a real basis of $\ggo '$
such that $\ad (z/c)=
\begin{pmatrix}\; 0 & 1 \\ -1 & 0
\end{pmatrix}$, hence $S11_{0,c} \cong \aff (\CC)$
for all $c>0$.

If $d=1$, then  $\ad (z)_{|\ggo '}$ has eigenvalues
 $1/2 \pm i \lambda /2, \;$ where $\lambda=
\sqrt{4c-1}$. Taking $z'=z/2$, there exists a real basis
of $\ggo '$ such that
$\ad (z')_{|\ggo '}$ takes the form $\begin{pmatrix} \; 1 & \lambda \\ -\lambda & 1
\end{pmatrix}$. Changing $z'$ to $z''=(z'-x)/ \lambda$, so that
$\ad (z'')= \begin{pmatrix}\; 0 & 1 \\ -1 & 0
\end{pmatrix}$, we conclude that $S11_{1,c} \cong \aff (\CC)$
for all $c$ such that $4c>1$.

\end{itemize}

\

Finally, in case the commutator ideal is three dimensional, we
establish the correspondence with Table 1 in \cite{O1}, p. 22.

\

\begin{center}
\small{
\begin{tabular}{|c|c|c|c|}\hline
$\begin{array}{c} \vspace{-.3cm} \\A1_{\lambda_1,\lambda_2} \, , \,\;\\ \vspace{-.3cm} \end{array} $
\scriptsize{$\lambda _1 \neq \lambda _2  \in \RR \backslash \{0,1\}$} &
$A1_{\lambda  , \overline{\lambda} } \, ,\;$  \scriptsize{$\text{Im}\, \lambda  \neq 0$} &
$A2_{\lambda} \, , \,\;$ \scriptsize{$\lambda \in \RR \backslash \{0,1\}$} &
$A3_{\lambda} \, , \, \;$ \scriptsize{$\lambda \in \RR \backslash \{0,1\}$} \\
\hline
$\begin{array}{c} \vspace{-.3cm} \\\mathfrak r _{4, \lambda _1 , \lambda _2}\\ \vspace{-.2cm} \end{array}$ &
$\mathfrak r'_{4, 1/ \text{Im}\,\lambda,\text{Re}\, \lambda  /
\text{Im}\, \lambda }$ & $\mathfrak r_{4,\lambda,\lambda}$ & $\mathfrak r _{4, \lambda}$ \\
 \hline
\end{tabular} }
\end{center}

\

\

\begin{center}
\small{
\begin{tabular}{|c| c|c|c|c|c|c|c|c|}\hline
 $\begin{array}{c} \vspace{-.3cm} \\ \;\;A4\;\; \\ \vspace{-.3cm} \end{array}$
 & $\;\;A5 \;\;$  & $\;\; A6 \;\; $ & $\;\; H1\;\;$ &
$H2$ & $H3$ & $H4$ &
$H5_{\lambda} \, , \,\;$ \scriptsize{$\lambda \in \RR \backslash \{0,1\}$} &
$H6 _{\lambda} \, , \,\;$ \scriptsize{$\lambda \in \CC \backslash \RR$}  \\
\hline
$\begin{array}{c} \vspace{-.3cm} \\ \mathfrak r _{4,1,1}\\ \vspace{-.2cm} \end{array}$
& $\mathfrak r _{4,1}$ & $\mathfrak r _{4}$ &  $ \mathfrak d _{4}$ &
$\mathfrak d' _{4, 0}$ & $\mathfrak d _{4, 1/2}$ & $\mathfrak h_4$ &
$\mathfrak d _{4, \lambda}$ & $\mathfrak d '_{4,-1/ \text{Im}\, \lambda}$ \\
\hline
\end{tabular} }
\end{center}

\

\section*{Appendix III - Some known results related to $4$-dimensional geometry} \label{known}

Using the characterization of homogeneous
manifolds of negative curvature given by Heintze in \cite{H} we
can conclude that
the following four dimensional Lie algebras do admit metrics with
negative curvature:
\begin{itemize}
\item $\mathfrak r _{4, \mu , \lambda}, \;\; 0< \mu \leq \lambda \leq 1$,
\item $\mathfrak r '_{4, \mu , \lambda}, \;\; \mu >0, \; \lambda > 0$,
\item $\mathfrak d_{4, \lambda}, \;\; 1/2 \leq \lambda <1$,
\item $\mathfrak d'_{4, \lambda}, \;\; \lambda >0$,
\item $\mathfrak h_4$.
\end{itemize}
Concerning non positive sectional curvature, we can mention
a result appearing in \cite{Dru}, where it is proved that a left invariant metric
with non positive curvature on a four dimensional solvable Lie group
either has geometric rank one or it comes from an inner product on
$\aff(\RR)\times \aff(\RR)$
 or $\RR \times \mathfrak r_{3,1}$, up to scaling.

We understand that the classification of rank one four dimensional homogeneous
spaces of non positive curvature is not known.
On the other hand, Jensen classified in \cite{J} the four dimensional
Lie algebras admitting Einstein metrics:
\begin{itemize}

\item $\RR \times \rr_{3,1}$,
\item $\mathfrak \rr_{4,1,1}$,
\item $\RR \times \mathfrak r '_{3,0}$,
\item $\mathfrak r_{4, \lambda , \lambda}', \;\;  \lambda > 0$.
\item $\mathfrak d_{4, \lambda}, \;\;  \lambda \geq 1/2$,

\end{itemize}

Among these, it follows from \cite{Al} that there are only two four dimensional Lie algebras
admitting Einstein metrics of non positive curvature: $\mathfrak r_{4,1,1}$ and $\mathfrak d_{4,1/2}$. Concerning left invariant anti-self-dual metrics on four dimensional Lie groups,
it was proved in \cite{D-S} (Theorem 1.6) that if a four dimensional Lie group admits such a metric, then its Lie algebra is one of the following:
\begin{itemize}
\item $\mathfrak d_{4, 1/2}$,
\item $\mathfrak d_{4, \lambda}', \; \lambda >0$.
\end{itemize}
It is proved in \cite{F2} that $\mathfrak d _{4,2}$ is the only four dimensional solvable Lie algebra admitting an almost K\"ahler structure
whose Ricci tensor is invariant with respect to the almost complex structure.

The classification of complex structures on four dimensional
solvable Lie algebras was carried out by Snow in \cite{SJ}, when
the commutator subalgebra has dimension one or two, and by Ovando
in \cite{O1}, when the commutator subalgebra is  three
dimensional. The classification of hypercomplex structures was
obtained in \cite{B1}.

Concerning the existence of symplectic structures, it is shown in
\cite{FG} that the solvmanifold obtained as a quotient of
$E(1,1)$, the simply connected Lie group  with Lie algebra
$\mathfrak e(1,1)$, by a lattice, admits a symplectic structure
but no complex structure. The classification of symplectic
structures on four dimensional Lie algebras is done in \cite{O2},
where the cohomology of all four dimensional solvable Lie algebras
is computed.

The hyper-K\"ahler metrics  conformal to  left invariant metrics
metric on  four dimensional Lie groups  were determined
 in \cite{B2}. It turns out that the solvable Lie groups appearing in this list
are those with Lie algebra $\RR^4, \; \aff (\CC), \; \rr_{4,1,1}$ or $\dd_{4,1/2}$.
It was proved in \cite{F1} that the cotangent bundle of a Lie group  with Lie algebra
$\aff(\CC)$ or $\rr_{4,1,1}$ also admits a metric conformal to a hyper-K\"ahler metric.

The determination of hypersymplectic structures on four dimensional Lie algebras was carried out in \cite{An}. According to this, the only Lie algebras admitting such a structure are $\RR^4,\,\RR\times\hh_3, \rr_{4,-1,-1}$ and $\dd_{4,2}$.


\

\

\

\begin{small} {\it E-mail addresses}: andrada@mate.uncor.edu,
barberis@mate.uncor.edu, idotti@mate.uncor.edu,

\hspace*{3cm}ovando@mate.uncor.edu \end{small}

\end{document}